\documentclass{article}
\usepackage{amssymb,amsmath,theorem,euscript}

\input{epsf}

\newcounter{sec}

\def\sm{\smallskip}

\newcounter{punct}[sec]

\def\punct{\refstepcounter{punct}{\arabic{sec}.\arabic{punct}.  }}

\def\COUNTERS{\addtocounter{sec}{1}
              \setcounter{punct}{0}
          \setcounter{equation}{0}
          \setcounter{theorem}{0}
                  }
\newtheorem{theorem}{Theorem}[sec]
\newtheorem{proposition}[theorem]{Proposition}
\newtheorem{lemma}[theorem]{Lemma}
\newtheorem{corollary}[theorem]{Corollary}
\newtheorem{conjecture}[theorem]{Conjecture}

\begin{document}

 \def\ov{\overline}
\def\wt{\widetilde}
\def\wh{\widehat}
 \newcommand{\rk}{\mathop {\mathrm {rk}}\nolimits}
\newcommand{\Aut}{\mathop {\mathrm {Aut}}\nolimits}
\newcommand{\Out}{\mathop {\mathrm {Out}}\nolimits}
\newcommand{\Abs}{\mathop {\mathrm {Abs}}\nolimits}
\renewcommand{\Re}{\mathop {\mathrm {Re}}\nolimits}
\renewcommand{\Im}{\mathop {\mathrm {Im}}\nolimits}
 \newcommand{\tr}{\mathop {\mathrm {tr}}\nolimits}
  \newcommand{\Hom}{\mathop {\mathrm {Hom}}\nolimits}
   \newcommand{\diag}{\mathop {\mathrm {diag}}\nolimits}
   \newcommand{\supp}{\mathop {\mathrm {supp}}\nolimits}
 \newcommand{\im}{\mathop {\mathrm {im}}\nolimits}
 \newcommand{\grad}{\mathop {\mathrm {grad}}\nolimits}
  \newcommand{\sgrad}{\mathop {\mathrm {sgrad}}\nolimits}
 \newcommand{\rot}{\mathop {\mathrm {rot}}\nolimits}
  \renewcommand{\div}{\mathop {\mathrm {div}}\nolimits}

\def\Br{\mathrm {Br}}
\def\Vir{\mathrm {Vir}}

 \def\Ham{\mathrm {Ham}}
\def\SL{\mathrm {SL}}
\def\Pol{\mathrm {Pol}}
\def\SU{\mathrm {SU}}
\def\GL{\mathrm {GL}}
\def\U{\mathrm U}
\def\OO{\mathrm O}
 \def\Sp{\mathrm {Sp}}
  \def\Ad{\mathrm {Ad}}
 \def\SO{\mathrm {SO}}
\def\SOS{\mathrm {SO}^*}
 \def\Diff{\mathrm{Diff}}
 \def\Vect{\mathfrak{Vect}}
\def\PGL{\mathrm {PGL}}
\def\PU{\mathrm {PU}}
\def\PSL{\mathrm {PSL}}
\def\Symp{\mathrm{Symp}}
\def\Cont{\mathrm{Cont}}
\def\End{\mathrm{End}}
\def\Mor{\mathrm{Mor}}
\def\Aut{\mathrm{Aut}}
 \def\PB{\mathrm{PB}}
\def\Fl{\mathrm {Fl}}
\def\Symm{\mathrm {Symm}} 
 \def\Herm{\mathrm {Herm}} 
  \def\SDiff{\mathrm {SDiff}} 
 
 \def\cA{\mathcal A}
\def\cB{\mathcal B}
\def\cC{\mathcal C}
\def\cD{\mathcal D}
\def\cE{\mathcal E}
\def\cF{\mathcal F}
\def\cG{\mathcal G}
\def\cH{\mathcal H}
\def\cJ{\mathcal J}
\def\cI{\mathcal I}
\def\cK{\mathcal K}
 \def\cL{\mathcal L}
\def\cM{\mathcal M}
\def\cN{\mathcal N}
 \def\cO{\mathcal O}
\def\cP{\mathcal P}
\def\cQ{\mathcal Q}
\def\cR{\mathcal R}
\def\cS{\mathcal S}
\def\cT{\mathcal T}
\def\cU{\mathcal U}
\def\cV{\mathcal V}
 \def\cW{\mathcal W}
\def\cX{\mathcal X}
 \def\cY{\mathcal Y}
 \def\cZ{\mathcal Z}
\def\0{{\ov 0}}
 \def\1{{\ov 1}}
 \def\frA{\mathfrak A}
 \def\frB{\mathfrak B}
\def\frC{\mathfrak C}
\def\frD{\mathfrak D}
\def\frE{\mathfrak E}
\def\frF{\mathfrak F}
\def\frG{\mathfrak G}
\def\frH{\mathfrak H}
\def\frI{\mathfrak I}
 \def\frJ{\mathfrak J}
 \def\frK{\mathfrak K}
 \def\frL{\mathfrak L}
\def\frM{\mathfrak M}
 \def\frN{\mathfrak N} \def\frO{\mathfrak O} \def\frP{\mathfrak P} \def\frQ{\mathfrak Q} \def\frR{\mathfrak R}
 \def\frS{\mathfrak S} \def\frT{\mathfrak T} \def\frU{\mathfrak U} \def\frV{\mathfrak V} \def\frW{\mathfrak W}
 \def\frX{\mathfrak X} \def\frY{\mathfrak Y} \def\frZ{\mathfrak Z} \def\fra{\mathfrak a} \def\frb{\mathfrak b}
 \def\frc{\mathfrak c} \def\frd{\mathfrak d} \def\fre{\mathfrak e} \def\frf{\mathfrak f} \def\frg{\mathfrak g}
 \def\frh{\mathfrak h} \def\fri{\mathfrak i} \def\frj{\mathfrak j} \def\frk{\mathfrak k} \def\frl{\mathfrak l}
 \def\frm{\mathfrak m} \def\frn{\mathfrak n} \def\fro{\mathfrak o} \def\frp{\mathfrak p} \def\frq{\mathfrak q}
 \def\frr{\mathfrak r} \def\frs{\mathfrak s} \def\frt{\mathfrak t} \def\fru{\mathfrak u} \def\frv{\mathfrak v}
 \def\frw{\mathfrak w} \def\frx{\mathfrak x} \def\fry{\mathfrak y} \def\frz{\mathfrak z} \def\frsp{\mathfrak{sp}}
 \def\bfa{\mathbf a} \def\bfb{\mathbf b} \def\bfc{\mathbf c} \def\bfd{\mathbf d} \def\bfe{\mathbf e} \def\bff{\mathbf f}
 \def\bfg{\mathbf g} \def\bfh{\mathbf h} \def\bfi{\mathbf i} \def\bfj{\mathbf j} \def\bfk{\mathbf k} \def\bfl{\mathbf l}
 \def\bfm{\mathbf m} \def\bfn{\mathbf n} \def\bfo{\mathbf o} \def\bfp{\mathbf p} \def\bfq{\mathbf q} \def\bfr{\mathbf r}
 \def\bfs{\mathbf s} \def\bft{\mathbf t} \def\bfu{\mathbf u} \def\bfv{\mathbf v} \def\bfw{\mathbf w} \def\bfx{\mathbf x}
 \def\bfy{\mathbf y} \def\bfz{\mathbf z} \def\bfA{\mathbf A} \def\bfB{\mathbf B} \def\bfC{\mathbf C} \def\bfD{\mathbf D}
 \def\bfE{\mathbf E} \def\bfF{\mathbf F} \def\bfG{\mathbf G} \def\bfH{\mathbf H} \def\bfI{\mathbf I} \def\bfJ{\mathbf J}
 \def\bfK{\mathbf K} \def\bfL{\mathbf L} \def\bfM{\mathbf M} \def\bfN{\mathbf N} \def\bfO{\mathbf O} \def\bfP{\mathbf P}
 \def\bfQ{\mathbf Q} \def\bfR{\mathbf R} \def\bfS{\mathbf S} \def\bfT{\mathbf T} \def\bfU{\mathbf U} \def\bfV{\mathbf V}
 \def\bfW{\mathbf W} \def\bfX{\mathbf X} \def\bfY{\mathbf Y} \def\bfZ{\mathbf Z} \def\bfw{\mathbf w}
 \def\R {{\mathbb R }} \def\C {{\mathbb C }} \def\Z{{\mathbb Z}} \def\H{{\mathbb H}} \def\K{{\mathbb K}}
 \def\N{{\mathbb N}} \def\Q{{\mathbb Q}} \def\A{{\mathbb A}} \def\T{\mathbb T} \def\P{\mathbb P} \def\G{\mathbb G}
 \def\bbA{\mathbb A} \def\bbB{\mathbb B} \def\bbD{\mathbb D} \def\bbE{\mathbb E} \def\bbF{\mathbb F} \def\bbG{\mathbb G}
 \def\bbI{\mathbb I} \def\bbJ{\mathbb J} \def\bbL{\mathbb L} \def\bbM{\mathbb M} \def\bbN{\mathbb N} \def\bbO{\mathbb O}
 \def\bbP{\mathbb P} \def\bbQ{\mathbb Q} \def\bbS{\mathbb S} \def\bbT{\mathbb T} \def\bbU{\mathbb U} \def\bbV{\mathbb V}
 \def\bbW{\mathbb W} \def\bbX{\mathbb X} \def\bbY{\mathbb Y} \def\kappa{\varkappa} \def\epsilon{\varepsilon}
 \def\phi{\varphi} \def\le{\leqslant} \def\ge{\geqslant}

\def\UU{\bbU}
\def\Mat{\mathrm{Mat}}
\def\tto{\rightrightarrows}

\def\F{\mathbf{F}}

\def\Gms{\mathrm {Gms}}
\def\Ams{\mathrm {Ams}}
\def\Isom{\mathrm {Isom}}

\def\Gr{\mathrm{Gr}}

\def\graph{\mathrm{graph}}

\def\O{\mathrm{O}}

\def\la{\langle}
\def\ra{\rangle}


 \def\ov{\overline}
\def\wt{\widetilde}

\renewcommand{\Re}{\mathop {\mathrm {Re}}\nolimits}
\def\Br{\mathrm {Br}}

 \def\Isom{\mathrm {Isom}}
 \def\Hier{\mathrm {Hier}}
\def\SL{\mathrm {SL}}
\def\SU{\mathrm {SU}}
\def\GL{\mathrm {GL}}
\def\U{\mathrm U}
\def\OO{\mathrm O}
 \def\Sp{\mathrm {Sp}}
  \def\GLO{\mathrm {GLO}}
 \def\SO{\mathrm {SO}}
\def\SOS{\mathrm {SO}^*}
 \def\Diff{\mathrm{Diff}}
 \def\Vect{\mathfrak{Vect}}
\def\PGL{\mathrm {PGL}}
\def\PU{\mathrm {PU}}
\def\PSL{\mathrm {PSL}}
\def\Symp{\mathrm{Symp}}
\def\ASymm{\mathrm{Asymm}}
\def\Asymm{\mathrm{Asymm}}
\def\Gal{\mathrm{Gal}}
\def\End{\mathrm{End}}
\def\Mor{\mathrm{Mor}}
\def\Aut{\mathrm{Aut}}
 \def\PB{\mathrm{PB}}
 \def\cA{\mathcal A}
\def\cB{\mathcal B}
\def\cC{\mathcal C}
\def\cD{\mathcal D}
\def\cE{\mathcal E}
\def\cF{\mathcal F}
\def\cG{\mathcal G}
\def\cH{\mathcal H}
\def\cJ{\mathcal J}
\def\cI{\mathcal I}
\def\cK{\mathcal K}
 \def\cL{\mathcal L}
\def\cM{\mathcal M}
\def\cN{\mathcal N}
 \def\cO{\mathcal O}
\def\cP{\mathcal P}
\def\cQ{\mathcal Q}
\def\cR{\mathcal R}
\def\cS{\mathcal S}
\def\cT{\mathcal T}
\def\cU{\mathcal U}
\def\cV{\mathcal V}
 \def\cW{\mathcal W}
\def\cX{\mathcal X}
 \def\cY{\mathcal Y}
 \def\cZ{\mathcal Z}
\def\0{{\ov 0}}
 \def\1{{\ov 1}}
 
 \def\frA{\mathfrak A}
 \def\frB{\mathfrak B}
\def\frC{\mathfrak C}
\def\frD{\mathfrak D}
\def\frE{\mathfrak E}
\def\frF{\mathfrak F}
\def\frG{\mathfrak G}
\def\frH{\mathfrak H}
\def\frI{\mathfrak I}
 \def\frJ{\mathfrak J}
 \def\frK{\mathfrak K}
 \def\frL{\mathfrak L}
\def\frM{\mathfrak M}
 \def\frN{\mathfrak N} \def\frO{\mathfrak O} \def\frP{\mathfrak P} \def\frQ{\mathfrak Q} \def\frR{\mathfrak R}
 \def\frS{\mathfrak S} \def\frT{\mathfrak T} \def\frU{\mathfrak U} \def\frV{\mathfrak V} \def\frW{\mathfrak W}
 \def\frX{\mathfrak X} \def\frY{\mathfrak Y} \def\frZ{\mathfrak Z} \def\fra{\mathfrak a} \def\frb{\mathfrak b}
 \def\frc{\mathfrak c} \def\frd{\mathfrak d} \def\fre{\mathfrak e} \def\frf{\mathfrak f} \def\frg{\mathfrak g}
 \def\frh{\mathfrak h} \def\fri{\mathfrak i} \def\frj{\mathfrak j} \def\frk{\mathfrak k} \def\frl{\mathfrak l}
 \def\frm{\mathfrak m} \def\frn{\mathfrak n} \def\fro{\mathfrak o} \def\frp{\mathfrak p} \def\frq{\mathfrak q}
 \def\frr{\mathfrak r} \def\frs{\mathfrak s} \def\frt{\mathfrak t} \def\fru{\mathfrak u} \def\frv{\mathfrak v}
 \def\frw{\mathfrak w} \def\frx{\mathfrak x} \def\fry{\mathfrak y} \def\frz{\mathfrak z} \def\frsp{\mathfrak{sp}}
 \def\bfa{\mathbf a} \def\bfb{\mathbf b} \def\bfc{\mathbf c} \def\bfd{\mathbf d} \def\bfe{\mathbf e} \def\bff{\mathbf f}
 \def\bfg{\mathbf g} \def\bfh{\mathbf h} \def\bfi{\mathbf i} \def\bfj{\mathbf j} \def\bfk{\mathbf k} \def\bfl{\mathbf l}
 \def\bfm{\mathbf m} \def\bfn{\mathbf n} \def\bfo{\mathbf o} \def\bfp{\mathbf p} \def\bfq{\mathbf q} \def\bfr{\mathbf r}
 \def\bfs{\mathbf s} \def\bft{\mathbf t} \def\bfu{\mathbf u} \def\bfv{\mathbf v} \def\bfw{\mathbf w} \def\bfx{\mathbf x}
 \def\bfy{\mathbf y} \def\bfz{\mathbf z} \def\bfA{\mathbf A} \def\bfB{\mathbf B} \def\bfC{\mathbf C} \def\bfD{\mathbf D}
 \def\bfE{\mathbf E} \def\bfF{\mathbf F} \def\bfG{\mathbf G} \def\bfH{\mathbf H} \def\bfI{\mathbf I} \def\bfJ{\mathbf J}
 \def\bfK{\mathbf K} \def\bfL{\mathbf L} \def\bfM{\mathbf M} \def\bfN{\mathbf N} \def\bfO{\mathbf O} \def\bfP{\mathbf P}
 \def\bfQ{\mathbf Q} \def\bfR{\mathbf R} \def\bfS{\mathbf S} \def\bfT{\mathbf T} \def\bfU{\mathbf U} \def\bfV{\mathbf V}
 \def\bfW{\mathbf W} \def\bfX{\mathbf X} \def\bfY{\mathbf Y} \def\bfZ{\mathbf Z} \def\bfw{\mathbf w}

 \def\R {{\mathbb R }} \def\C {{\mathbb C }} \def\Z{{\mathbb Z}} \def\H{{\mathbb H}} \def\K{{\mathbb K}}
 \def\N{{\mathbb N}} \def\Q{{\mathbb Q}} \def\A{{\mathbb A}} \def\T{\mathbb T} \def\P{\mathbb P} \def\G{\mathbb G}
 \def\bbA{\mathbb A} \def\bbB{\mathbb B} \def\bbD{\mathbb D} \def\bbE{\mathbb E} \def\bbF{\mathbb F} \def\bbG{\mathbb G}
 \def\bbI{\mathbb I} \def\bbJ{\mathbb J} \def\bbL{\mathbb L} \def\bbM{\mathbb M} \def\bbN{\mathbb N} \def\bbO{\mathbb O}
 \def\bbP{\mathbb P} \def\bbQ{\mathbb Q} \def\bbS{\mathbb S} \def\bbT{\mathbb T} \def\bbU{\mathbb U} \def\bbV{\mathbb V}
 \def\bbW{\mathbb W} \def\bbX{\mathbb X} \def\bbY{\mathbb Y} \def\kappa{\varkappa} \def\epsilon{\varepsilon}
 \def\phi{\varphi} \def\le{\leqslant} \def\ge{\geqslant}

\def\UU{\bbU}
\def\Mat{\mathrm{Mat}}
\def\tto{\rightrightarrows}

\def\Gr{\mathrm{Gr}}

\def\B{\bfB} 

\def\graph{\mathrm{graph}}

\def\O{\mathrm{O}}

\def\gl{\mathfrak{gl}}

\def\la{\langle}
\def\ra{\rangle}

\begin{center}
\Large\bf
Restriction of representations of $\GL(n+1,\C)$ to $\GL(n,\C)$ and  action
of the Lie overalgebra	

\bigskip

\sc Yury A. Neretin
\footnote{Supported by the grant FWF, P25142, P28421}
\end{center}

\section{The statement}

\COUNTERS

{\small{\sc Abstract.} Consider a restriction of an irreducible
finite dimensional holomorphic representation of $\GL(n+1,\C)$ to the subgroup
$\GL(n,\C)$. We write explicitly formulas for generators of the Lie algebra
$\frg\frl(n+1)$ in the direct sum of representations of $\GL(n,\C)$.
Nontrivial generators act as differential-difference operators, the differential part
has order $n-1$, the difference part acts on the space of parameters (highest weights)
of representations. We also formulate a conjecture about unitary principal series of 
$\GL(n,\C)$.}

\bigskip

{\bf\punct The Gelfand--Tsetlin formulas.}
It is well known that restrictions of finite dimensional holomorphic representations
of the general linear group $\GL(n,\C)$ to the subgroup $\GL(n-1,\C)$ is multiplicity
free. Considering a chain of restrictions
$$\GL(n,\C)\supset\GL(n-1,\C)\supset \GL(n-2,\C)\supset\dots\supset \GL(1,\C)$$
we get a canonical decomposition of 
our representation into a direct sum of one-dimensional subspaces.
Taking a vector in each line we get a basis of
the representation. In \cite{GTs1} Gelfand and Tsetlin announced  formulas for  action
of generators of the Lie algebra $\frg\frl(n)$ in this basis. It turns out that
the Lie algebra $\frg\frl(n)$ acts by difference operators in the space of functions on a certain  convex polygon in 
the lattice $\Z^{n(n-1)/2}$. In particular, this gives an explicit realization of representations
of the Lie algebra $\gl(n)$.

For various proofs of the Gelfand--Tsetlin formulas and
constructions of the bases, see \cite{BB}, \cite{Zhe-book},
\cite{AST}, \cite{Mick}, \cite{Zhe-red}, \cite{NT}, \cite{Molev2}, \cite{Molev3} (see more references in \cite{Molev2}), 
for other classical groups, see
\cite{GTs2},   \cite{Molev2}, \cite{Sht}, for applications to infinite-dimensional representations,
see \cite{Olsh}, \cite{Mol2}, \cite{Gra}, for formulas on the group
level, see \cite{GG}, \cite{Gra}. There are  many other continuations
of this story.
However, our standpoint is slightly different.

\sm

{\bf \punct Actions of overalgebras in the spectral decompositions.}
We can write
images of many operators under the classical  Fourier transform. It is commonly accepted
that
Plancherel decompositions of representations are higher analogs of 
the Fourier transform. 

Consider a group 
$G$ and its subgroup $H$. Restrict an irreducible unitary representation of 
$G$ to $H$. Generally, an explicit spectral decomposition of
the restriction seems hopeless problem. However, there is a collection of explicitly solvable problems
of this type%
\footnote{Also, some explicitly solvable spectral problems 
in representation theory can be regarded as special cases of the restriction problem.
In particular, a decomposition of $L^2$ on a classical pseudo-Riemannian symmetric space $G/H$ can be regarded
as a special case of the restriction of a Stein type principal series of a certain overgroup $\wt G$ to
the symmetric subgroup $G$, \cite{Ner-Stein}. So, for $L^2$ on symmetric spaces the problem 
of action of the overalgebra discussed below
makes sense.}.
In \cite{Ner-over} it was conjectured that in such cases the action
of the Lie algebra $\frg$ can be written explicitly as differential-difference
operators. In fact, in \cite{Ner-over} there was considered the tensor product
$V_s\otimes V_s^*$
of a highest weight and a lowest weight unitary representations 
of $\SL(2,\R)$ (i.e., $G\simeq \SL(2,\R)\times \SL(2,\R)$
and $H\simeq \SL(2,\R)$ is the diagonal). This representation
is a multiplicity free integral over the principal series of $\SL(2,\R)$.
It appears that the action of  all generators of the Lie algebra $\frs\frl(2,\R)\oplus \frs\frl(2,\R)$
 in the spectral decomposition can be written explicitly
 as differential-difference operators. The differential part of these
 operators has order two and the difference operators are difference operators 
 in the imaginary direction%
 \footnote{On Sturm--Liouville in imaginary direction, see \cite{Ner-imaginary} and further references in that paper,
 see also \cite{Gro}.
 The most of known appearances of such operators are related to representation theory and spectral
 decompositions of unitary representations.}.
 
 Molchanov \cite{Mol1}--\cite{Mol4} solved several problems of this type 
 related to rank-one symmetric spaces%
 \footnote{In particular, he examined the restrictions of maximally degenerate principal
 series for the cases $\GL(n+1,\R)\supset \OO(n,1)$, $\OO(p,q)\supset \OO(p-1,q)$.}.
 In all his cases the Lie overalgebra
 acts by differential-difference operators; difference operators act in the imaginary dimension.
 However,  the order of the differential part in some cases is 4. In \cite{Ner-jfaa}, \cite{Ner-operational}
 the author considered the Fourier transform on the group $\GL(2,\R)$. There were evaluated 
 the operators in the target space corresponding to multiplications by matrix elements and 
 partial derivatives with respect to matrix elements. Formulas have similar structure, but in \cite{Ner-operational}
 they are simpler than in previous cases. 
 
 In the present  paper, we consider  restrictions of holomorphic finite-dimensional representations
 $\GL(n+1,\C)$ to $\GL(n,\C)$ and write explicit formulas 
 for the action of the overalgebra in the spectral decomposition.
 Also, we formulate a conjecture concerning restrictions of unitary principal series.

 \sm

{\bf\punct Notation. The group $\GL(n,\C)$ and its subgroups.}
Denote by $\GL(n,\C)$ the group of invertible complex matrices of size
$n$.  By $\gl(n)$ we denote its Lie algebra, i.e. the Lie algebra
of all matrices of size $n$. Let $\U(n)\subset \GL(n,\C)$ be the subgroup
of unitary matrices.

By $E_{ij}$ we denote  the matrix whose $ij$'s entry is 1 and other entries are 0.
The unit matrix is denoted by $1$ or $1_n$, i.e. $1=\sum_j E_{jj}$.
Matrices $E_{ij}$ can be regarded as generators
of the Lie algebra $\gl_n$.

Denote by $N^+=N_n^+\subset \GL(n,\C)$ the subgroup of all strictly upper triangular matrices of
size $n$, i.e., matrices of the form
$$
Z=\begin{pmatrix}
1&z_{12}& \dots& z_{1(n-1)}&z_{1n}\\
0&1&\dots & z_{2(n-1)}&z_{2n}\\
\vdots&\vdots &\ddots &\vdots&\vdots\\
0&0&\dots&1&z_{(n-1)n}\\
0&0&\dots&0&1
\end{pmatrix}.
$$
By $N_n^-$ we denote the group of strictly lower triangular matrices. 

By $B_n^-\subset \GL(n,\C)$ we denote the subgroup
of lower triangular matrices, $g\in B_n^-$ if $g_{ij}=0$ if
$i<j$. Let $\Delta$ be the subgroup of diagonal matrices,
we denote its elements as $\delta=\diag(\delta_1, \dots, \delta_n)$. 

\sm

Next, we need a notation for sub-matrices. For a matrix
$X$  denote by 
$$[X]_{\alpha\beta}$$
 the left 
upper corner of $X$ of size $\alpha\times \beta$.
Let $I=\{i_1,\dots, i_\alpha\}$, $J=\{j_1, \dots, j_\beta\}$ be collections of integers,
we assume that their elements are ordered as $0<i_1<i_2<\dots i_\alpha$, $J:1\le j_1\le \dots\le j_\beta$.
Denote by 
\begin{equation} 
\left[X\begin{pmatrix} I\\J \end{pmatrix}  \right]
=
\left[X\begin{pmatrix}
i_1,\dots, i_\alpha\\j_1,\dots, j_\beta 
\end{pmatrix}  \right]
\label{eq:submatrix}
\end{equation} 
the matrix composed of entries
$x_{i_\mu,j_\nu}$.
By 
\begin{equation}
[X(I)]=\bigl[X(i_i,\dots,i_\alpha)\bigr]
\label{eq:submatrix-1}
\end{equation}
we denote matrix composed of  $i_i$,\dots, $i_\alpha$-th rows of $X$
(the order is not necessary increasing, also we allow coinciding rows).

\sm 

{\bf \punct Holomorphic representations of $\GL(n,\C)$.}
Recall that irreducible finite dimensional holomorphic representations
$\rho$ of the group $\GL(n,\C)$ are enumerated by collections of integers
({\it signatures})
$$\bfp:\,p_1\ge \dots \ge p_n.$$
This means that there is a cyclic vector $v$ (a {\it highest weight vector}) such that 
\begin{align*} 
\rho(Z)v&=v\qquad \text{for each
$Z\in N_n^+$},\\
\rho(\delta) v&=\prod_{j=1}^n \delta_j^{p_j}v, \qquad\text{where $\delta=\diag(\delta_1, \dots, \delta_n)$} 
.\end{align*}
Denote such representation by $\rho_\bfp=\rho^n_\bfp$.
Denote the set of all signatures  by $\Lambda_n$.

The dual (contragradient) representation to $\rho_\bfp$
has the signature 
\begin{equation}
\bfp^*:=(-p_n,\dots,-p_1).
\label{eq:dual-signature}
\end{equation}

\sm

{\bf\punct Realizations of representations by differential operators.}  
Recall a model for irreducible finite dimensional holomorphic representations
for details, see Section 2).
We consider the space $\Pol(N_n^+)$ of
polynomials in the variables $z_{ij}$. The generators of
the Lie algebra $\gl(n)$ act 
in this space via differential operators
\begin{align}
E_{kk}&=-\sum_{i<k} z_{ik}\partial_{ik}+p_k  +\sum_{j>k} z_{kj}\partial_{kj},
\label{eq:gener1}
\\
E_{k(k+1)}&=\partial_{k(k+1)}+\sum_{i<k} z_{ik}\partial_{i(k+1)},
\label{eq:gener2}
\\
E_{(k+1)k}&=\sum_{i<k} z_{i(k+1)}\partial_{ik}+(p_{k+1}-p_k)z_{k(k+1)}-
z_{k(k+1)} \sum_{j>k} z_{kj} \partial_{kj}+
\nonumber 
\\&\qquad\qquad\qquad +\sum_{m>k+1}
\det\begin{pmatrix}
z_{k(k+1)} &z_{km}\\
1&z_{(k+1)m}
\end{pmatrix}
\partial_{(k+1)m},
\label{eq:gener3}
\end{align}
where 
$$
\partial_{kl}:=\frac{\partial}{\partial z_{kl}}.
$$
All other generators can be expressed in the terms of $E_{kk}$, $E_{k(k+1)}$,
$E_{(k+1)k}$ (also it is easy to write explicit formulas as it is explained in Section 2). In this way, we get a representation of the Lie algebra
$\gl(n)$.

 There exists a unique finite dimensional subspace $V_\bfp$ invariant
with respect to such operators. 
The representation $\rho_\bfp$ is realized in this space.
 The
highest weight vector is $f(Z)=1$. In the next subsection we describe this space more explicitly.

\sm

{\sc Remark.} This approach arises in \cite{GN}.
Of course,  this  model is a coordinatization of the constructions of representations of $\GL(n,\C)$
in sections of line bundles over the flag space $B_n^-\setminus \GL(n,\C)$
(as in the Borel--Weil--Bott theorem).
The space $N_n^+\simeq \C^{n(n-1)/2}$ is an open dense chart on this space,
the terms with  $\partial$ in (\ref{eq:gener1})--(\ref{eq:gener3})
correspond to vector fields on the flag space,  the zero order 
terms are corrections corresponding to the action in the bundle. 
Elements of the space $V_\bfp^n$ are precisely polynomials, which are holomorphic 
as section of bundle on the whole flag space.
However, we need explicit formulas and prefer
a purely coordinate language.
\hfill $\boxtimes$

\sm

{\bf\punct Descriptions of the space $V_\bfp$. Zhelobenko operators.}

\sm

a) {\sc Description-1.}
Denote by $d\dot Z$ the standard Lebesgue measure on $N_n^+$
$$
d\dot Z=\prod_{k<l} d\Re z_{kl}\,d\Im z_{kl}.
$$
Denote by $d\mu(Z)$ the 
measure on $N_n^+$ given by the formula 
\begin{equation}
d\mu_\bfp(Z)=d\mu^n_\bfp(Z)=\prod_{j=1}^{n-1} \det\bigl([Z]_{jn}([Z]_{jn})^*\bigr)^{-(p_j-p_{j-1})-2}\,d\dot Z
.
\label{eq:mu}
\end{equation}

\begin{proposition}
	\label{prop:invariance-of-measure}
	{\rm a)}
 $V_\bfp^n=L^2(N_n^+,\mu^n_\bfp) \cap \Pol(N_n^+)$.	

\sm

{\rm b)} The $L^2$-inner product in $V^n_\bfp$ is $\U(n)$-invariant. 
\end{proposition}

Proof is given in Section 2.

\sm

{\sc Description-2.} We define the {\it Zhelobenko operators} by
\begin{equation}
R_{km}:=\partial_{km}+\sum_{j>m} z_{mj}\partial_{kj}, \qquad m>k
.
\label{eq:zhelobenko}
\end{equation}

\begin{theorem}
\label{th:zhelobenko}
The space $V_\bfp^n$ consists of polynomials satisfying the conditions
$$
\begin{cases}
R_{j(j+1)}^{p_j-p_{j+1}+1} f(Z)=0.
\end{cases}
$$ 	
	\end{theorem}

See Zhelobenko, \cite{Zhe-book}, \cite{Zhe-paper}, \S6, Theorem 2, \cite{Zhe-principal}, Theorem 48.7.

\sm

 {\sc Description-3.} There is one more description of the space
 $V_\bfp^n$, which can be used with coordinate language.

\begin{proposition}
	\label{prop:reproducing}
	The space $V_\bfp^n$ coincides with the space determined by the reproducing kernel
	$$
	\prod_{j=1}^{n-1} \det\left([Z]_{jn}[Z]_{jn}^* \right)^{p_j-p_{j-1}}
	.
	$$
		\end{proposition}
		
On Hilbert spaces determined by reproducing kernels, see, e.g., \cite{Ner-gauss}, Sect. 7.1.	
The proposition is more-or-equivalent to the Borel--Weil theorem.

\sm

{\bf\punct The restriction of representations $\GL(n+1,\C)$ to $\GL(n,\C)$.}
Consider the representation $\rho^{n+1}_\bfr$ of $\GL(n+1,\C)$ with a signature
$\bfr=(r_1,\dots,r_{n+1})$. It is well-known (see Gelfand, Tsetlin \cite{GTs1}) that the restriction
of $\rho^{n+1}_\bfr$
to $\GL(n,\C)$ is multiplicity free and is a direct sum of all representations
$\rho_\bfq$ of $\GL(n,\C)$ with signatures satisfying the {\it interlacing
conditions}
\begin{equation}
r_1\ge q_1\ge r_2\ge q_2\ge r_3\ge \dots\ge q_n\ge r_{n+1}.
\label{eq:interlacing}
\end{equation}

Our purpose is to write explicitly the action of Lie algebra
$\gl(n+1)$ in the space 
$\oplus V_\bfq^n$.

\sm

{\bf \punct Normalization.%
	\label{ss:normalization}}
First,  we intend to write a $\GL(n,\C)$-invariant pairing
$$J_{\bfp,\bfq}:V^{n+1}_\bfp\times V^n_\bfq\to \C.$$
as an integral
\begin{equation}
J_{\bfp,\bfq} (\phi,\psi)=\int_{N_{n+1}^+\times N_n^+} \ov{ L_{\bfp,\bfq} (U,Z)} \phi(Z)\,\psi(U)  \,d\mu_\bfp^{n+1} (Z)
\,d\mu^n_\bfq(U),
\label{eq:Lpq}
\end{equation}
where $Z$ ranges in $N_{n+1}^+$, $U$ ranges in $N_n^+$, and the kernel
$$L_{\bfp,\bfq} (U,Z)
\in  V^{n+1}_\bfp\otimes V^n_\bfq
$$
is a polynomial in the variables
$z_{ij}$, $ u_{kl}$.

Denote  $Z^{cut}:=[Z]_{(n+1)n}$, this matrix is obtained from $Z$ by cutting of
the last column.  Denote by $U^{ext}$ the $n\times(n+1)$-matrix obtained 
from $U$ by adding the zero last column,
\begin{equation}
U^{ext}=\begin{pmatrix}
1&u_{12}&\dots & u_{1n}&0\\
0&1&\dots&u_{2n}&0\\
\vdots&\vdots&\ddots &\vdots\\
0&0&\dots&1&0
\end{pmatrix}.
\label{U-ext}
\end{equation}
Consider the $(n+1)\times (n+1)$-matrix
composed from the first  $(n+1-\alpha)$ rows of the matrix $Z$ and the first $\alpha$
rows of the matrix $U^{ext}$.
Denote by $\Phi_\alpha$ the determinant of this matrix:
\begin{equation}
\Phi_\alpha(Z,U)=\det
\begin{pmatrix}
[Z]_{ (n+1-\alpha) (n+1)}
\\
[U^{ext}]_{\alpha(n+1)}
\end{pmatrix},
\qquad
\end{equation}

 Consider
the $n\times n$-matrix composed of the first $n-\alpha$ rows of the matrix
$Z^{cut}$ and the first  $\alpha$ rows of the matrix $U$,
denote by $\Psi_{\alpha}$ its determinant:
\begin{equation}
\Psi_\alpha(Z,U)=
\det
\begin{pmatrix}
[Z^{cut}]_{ (n-\alpha) n}\\
[U]_{\alpha n}
\end{pmatrix}.
\label{eq:Phi-Psi}
\end{equation}

\begin{proposition}
	\label{prop:kernel}
	Consider signatures $\bfp=(p_1,\dots,p_{n+1})$, $\bfq=(q_1,\dots, q_n)$
	such that $\bfp$ and $\bfq^*$ are interlacing. Then
	the expression  {\rm(\ref{eq:Lpq})} with the kernel 
	\begin{equation}
	L_{\bfp,\bfq}(U,Z)=
	\Phi_1^{p_n+q_1}\Psi_1^{-p_{n}-q_2} \Phi_2^{p_{n-1}+q_2} 
	\Psi_2^{-p_{n-1}-q_3}
	\dots \,\Phi_n^{p_1+q_n}.
	\label{eq:kernel}
	\end{equation} 	
	determines a $\GL(n,\C)$-invariant nonzero pairing between $V_\bfp^{n+1}$ and  $V_\bfq^n$.
\end{proposition}

For instance, for $n=3$ we get
{\small
\begin{multline*}
L_{\bfp,\bfq}(U,Z)=
\det\begin{pmatrix}
1& z_{12}& z_{13}& z_{14}\\
0& 1& z_{23}& z_{24}\\
0&0&1& z_{34}\\
1&u_{12}&u_{13}&0
\end{pmatrix}^{p_3+q_1}
\det\begin{pmatrix}
1& z_{12}& z_{13}\\
0&1& z_{23}\\
1&u_{12}&u_{13}\\
\end{pmatrix}^{-p_3-q_2}
\times\\
\det\begin{pmatrix}
1& z_{12}& z_{13}& z_{14}\\
0& 1& z_{23}& z_{24}\\
1&u_{12}&u_{13}&0\\
0&1&u_{23}&0
\end{pmatrix}^{p_2+q_2}
\!\!\!\!\!\!
\det\begin{pmatrix}
1& z_{12}& z_{13}\\
1&u_{12}&u_{13}\\
0&1&u_{23}\\
\end{pmatrix}^{-p_2-q_3}
\!\!\!\!\!\!
\det\begin{pmatrix}
1& z_{12}& z_{13}& z_{14}\\
1&u_{12}&u_{13}&0\\
0&1&u_{23}&0\\
0&0&1&0
\end{pmatrix}^{p_1+q_3}
.
\end{multline*}} 

Next we pass  to the dual signature 
$$\boxed{\bfr=\bfp^*}$$
(below $\bfp$ and $\bfq$ are rigidly linked by this restraint),
and represent $L_{\bfp,\bfq}$ in the form
	\begin{equation}
	L^{\bfr}_\bfq(U,Z):=L_{\bfp,\bfq}(U,Z)=
	\Phi_1^{q_1-r_2}\Psi_1^{r_2-q_2} \Phi_2^{q_2-r_3} 
	\Psi_2^{r_3-q_3}
	\dots \Phi_n^{q_n-r_{n+1}}
	\label{eq:kernel-1}
	.
	\end{equation}

	\sm

{\bf \punct Action of the overalgebra.%
\label{ss:ao}}  
Fix a signature $\bfr=\bfp^*\in \Lambda_{n+1}$ as in the previous subsection.
Consider the space
\begin{equation}
\frV_{\bfr}:=\bigoplus_{q_1, \dots, q_n:\,r_j\ge q_j\ge r_{j-1} } V_\bfq^n.
\label{eq:frV}
\end{equation}
We can regard elements of this space as 
 'expressions'
$$
F(U,\bfq)
$$
 of $n(n-1)/2$ continuous   variables $u_{kl}$, where $1\le k<l\le n$, and of
integer variables
$q_1$, \dots, $q_n$. Integer variables range in
the domain $r_j\ge q_j\ge r_{j+1}$. More precisely, for a fixed $\bfq$
the expression $F(U,\bfq)$  is a polynomial in the variables $u_{kl}$,
moreover, this polynomial  is contained in $V_\bfq$.

Our family of forms $L_{\bfp,\bfq}$ determines a duality between
$V_{\bfr^*}$ and  the space $\frV_{\bfr}$, hence we get a canonical  identification of
$V_\bfr$ and $\frV_{\bfr}$.
Therefore, we have a canonically defined action
of the Lie algebra $\frg\frl(n+1)$ in $\frV_\bfr$.
We preserve the notation $E_{kl}$ for operators in  $V_\bfp$ 
and denote operators in $\frV_\bfr$ by $F_{kl}$.
 For $1\le k,l\le n$ the operators $F_{kl}^n$ act
 in the space $\frV_\bfr$ by the first order differential
 operators  in $U$ with coefficients depending on $\bfq$
 according the standard formulas (see (\ref{eq:gener1})--(\ref{eq:gener3})). For instance,
\begin{align}
F_{kk}&=-\sum_{i<k} u_{ik}\partial_{ik}+q_k  +\sum_{j>k} u_{kj}\partial_{kj},
\label{eq:eu-1}
\\
F_{k(k+1)}&=\partial_{k(k+1)}+\sum_{i<k} u_{ik}\partial_{i(k+1)},
\\
F_{(k+1)k}&=\sum_{i<k} u_{i(k+1)}\partial_{ik}+(q_{k+1}-q_k)u_{k(k+1)}-
u_{k(k+1)} \sum_{j>k} u_{kj} \partial_{kj}+
\\&\qquad\qquad\qquad +\sum_{m>k+1}
\det\begin{pmatrix}
u_{k(k+1)} &u_{km}\\
1&u_{(k+1)m}
\end{pmatrix}
\partial_{(k+1)m},
\label{eq:eu-3}
\end{align}
where $\partial_{kl}:=\frac{\partial}{\partial u_{kl}}$.

The purpose of this work is  to present formulas for the generators
 $F_{1(n+1)}$, $F_{(n+1)n}$. Together with (\ref{eq:eu-1})--(\ref{eq:eu-3})
 they generate the Lie algebra $\frg\frl(n+1)$, formulas for the remaining 
 generators
 $F_{j(n+1)}$, $F_{(n+1)j}$ consist of similar aggregates as below, but are longer.

Denote by $T_j^\pm$ the following difference operators
$$
T_j^\pm F(\dots,q_j, \dots)= F(\dots,q_j\pm1, \dots)
$$
(the remaining variables do not change). We will write expressions,
which are polynomial in $u_{kl}$, $\partial_{kl}$,
linear in $T_j^\pm$ and rational in $q_j$. These expressions  satisfy
the commutation relations in $\gl(n)$ and preserve the space $\frV_\bfr$.

For $1\le k <l\le n$ denote
by $[k,l]$ the set
$$
[k,l]:=\{k,k+1,\dots,l\}.
$$ 
For a set
$$
I: 1\le i_1<i_2<\dots<i_m
$$ 
we write
\begin{equation}
I\triangleleft [k,l] \qquad \text{if $i_1=k$, $i_m=l$.}
\label{eq:triangle}
\end{equation}
More generally, if $I\subset J\subset[1,n]$,
we write $I\triangleleft J$ if the minimal (resp. maximal) 
element of $I$ coincides with the minimal (resp. maximal) 
element of $J$.

For any $I \triangleleft [\alpha,\beta]$  we define the operator
$$
R_I:=R_{i_1i_2}R_{i_2i_3}\dots R_{i_{m-1}i_m}= R_{\alpha i_2}R_{i_2i_3}\dots R_{i_{m-1}\beta},
$$
where $R_{ij}$ are the Zhelobenko operators.
We also set
$$
R_{kk}:=1.
$$

\begin{theorem}
{\rm	a)} The generator $F_{1(n+1)}$ acts  by the formula
	\begin{equation}
	F_{1(n+1)}=\sum_{m=1}^n  
	A_m
	\Bigl(\sum_{I\triangleleft [1,m]}\,\,
	 \prod_{l\in [1,m]\setminus I} (q_m-q_l+l-m) \cdot R_I
	\Bigr) T^-_m,
		\label{eq:formula-1}
	\end{equation}
	where coefficients $A_m$ are given by
	$$
	A_m=
	\frac{\prod_{j=m+1}^{n+1} (q_m-r_j+j-m-1)}
	{\prod_{\alpha\ne m} (q_m-q_\alpha+\alpha-m)}
	.$$

	{\rm	b)} The generator $F_{(n+1)n}$ is given by the formula
	\begin{equation}
	F_{(n+1)n}=-
	\sum_{m=1}^n B_m \Bigl(\sum_{I\triangleleft [m,n]}\,\,
	 \prod_{l\in [m,n]\setminus I} (q_m-q_l-m+l+1) \cdot R_I \Bigr)T_m^+	
	 ,
	 \label{eq:formula-2}
	\end{equation}
where
$$
B_m=
\frac{\prod_{j=1}^{m}(q_m-r_{j}+j-m)}{\prod_{\alpha\ne m}(q_m-q_\alpha-m+\alpha)}.
$$	
\end{theorem}

{\bf \punct Further structure of the paper.}
Section 2 contains preliminaries on holomorphic  representations of $\GL(n,\C)$.
In Section 3, we verify the formula for the kernel $L(Z,U)$.
Our main  statement is equivalent to a verification of  differential-difference equations
for the kernel $L(Z,U)$, this is done in Section 4. In Section 5, we formulate
a conjecture about analog of our statement for unitary representations.

\sm

{\bf Acknowledgements.} Fifteen years ago the topic of the paper
was one of aims of a joint project with M.~I.~Graev,
which was not realized in that time (I would like to emphasis his nice paper \cite{Gra}).
I am grateful to him and also to V.~F.~Molchanov
for discussions of the problem.

\section{Holomorphic representations of $\GL(n,\C)$.}

\COUNTERS

{\bf\punct Realization in the space of functions on $\GL(n,\C)$.%
\label{ss:realization-1}} For details and proof, see
\cite{Zhe-book}.
We say that a function $f(g)$ on $\GL(n,\C)$ is a polynomial function if it can be expressed
as a polynomial expression in matrix elements $g_{ij}$ and $\det(g)^{-1}$.
Denote the space of polynomial functions by
$\C[\GL(n,\C)]$. The group $\GL(n,\C)\times \GL(n,\C)$ acts in $\C[\GL(n,\C)]$
by the left and right shifts`:
$$
\lambda_{left-right}(h_1,h_2) f(g)= f(h_1^{-1}g h_2).
$$
This representation is a direct sum 
$$
\lambda \simeq \oplus_{\bfp\in\Lambda_n} \,\rho_{\bfp^*}\otimes \rho_\bfp.
$$
A summand $V_{\bfp^*}\otimes V_{\bfp}$ is the $\GL(n,\C)\times \GL(n,\C)$-cyclic span
of the vector
$$
\Delta[\bfp]=\prod_{j=1}^{n-1} \det[g]_{jj}^{p_j-p_{j+1}} \cdot (\det g)^{p_n}
.$$
If $h_1\in N_n^-$ and $h_2\in N_n^+$, then  
$\lambda(h_1,h_2)\Delta[\bfp]=\Delta[\bfp]$.

\sm

Next, consider the space  of polynomial functions invariant with respect to left shifts 
on elements of $N_n^-$
\begin{equation}
f(h^{-1}g)=f(g).
\label{eq:left-invariance}
\end{equation}
Equivalently, we consider the space of polynomial functions 
on $\C[N_n^-\setminus \GL(n,\C)]$.
Each $\rho_{\bfp^*}$ has a unique $N_n^-$-invariant vector, therefore
$\C[N_n^-\setminus \GL(n,\C)]$ is a multiplicity free direct sum
$\oplus_\bfp \rho_\bfp$. 

Fix a signature $\bfp$ and consider the space of $N_n^-$-invariant 
functions $H_\bfp$ such that for any diagonal matrix
$\delta$, we have
$$
f(\delta g)=f(g)\prod_{j=1}^n \delta_j^{p_j}.
$$
The group $\GL(n,\C)$ acts in $H_\bfp$ by the right shifts,
\begin{equation}
\lambda_{right}(h) f(g)=f(gh)
\label{eq:right-regular}
.\end{equation}
This representation is irreducible and equivalent to
$\rho_\bfp$. The vector $\Delta_\bfp$ is its highest weight vector.

\sm

{\bf \punct Realization in the space of functions on $N_n^+$.%
\label{ss:realization-2}}
 An element $g\in\GL(n,\C)$ 
satisfying the condition
$$
\det [g]_{jj}\ne 0\qquad \text{for all $j$}
$$
admits a unique Gauss decomposition 
$$
g=bZ,\qquad\text{where $b\in B_n^-$, $Z\in N_n^+$}
.$$
Notice that
\begin{align}
b_{jj}&=\frac {\det [g]_{jj}}{\det [g]_{(j-1)(j-1)}}
,
\label{eq:b-Gauss}
\\
z_{ij}&=\frac{\det\left[g\begin{pmatrix}1&2&\dots&i-1&j\\1&2&\dots&i&j \end{pmatrix}\right]}
{\det\left[g\begin{pmatrix}1&2&\dots&i&j\\1&2&\dots&i&j \end{pmatrix}\right]}
.
\label{eq:z-gauss}
\end{align}

Restrict a function $f\in H_\bfp$ to the subgroup $N_n^+$. We get a polynomial 
in the variables $z_{ij}$, where $i<j$. By the Gauss decomposition, $f$
is uniquely determined by this restriction. Therefore the space $V_\bfp$
can be regarded as a subspace of the space of polynomial in $z_{ij}$.
The description of this space is given by the Zhelobenko Theorem
\ref{th:zhelobenko}.

Denote the factors in the Gauss decomposition of $Zg$ in the following way
\begin{equation}
Zg=:b(Z,g)\cdot Z^{[g]}, \qquad \text{where $b(Z,g)\in B_n^-$, $Z^{[g]}\in N_n^+$},
\label{eq:urodina}
\end{equation}
Then
\begin{equation}
\rho_\bfp (g)\,f(Z)=f(Z^{[g]})  \,\prod_{j=1}^{n} b_{jj}(Z,g)^{p_j}
,
\label{eq:group-action}
\end{equation}
where $b_{jj}(Z,g)$ are the diagonal matrix elements of the lower triangular matrix
$b(Z,g)$. 

We present formulas for transformations $Z\mapsto Z^{[g(t)]}$,
where $g(t)$ are standard one-parametric subgroups in $\GL(n,\C)$.

\sm

1) Let $g(t)=\exp (1+t E_{kl})$, where $k<l$. Then
we have transformation
\begin{align*}
z_{kl}&\mapsto z_{kl}+t;\\
z_{ml}&\mapsto z_{ml}+t z_{kl}, \qquad m<k;\\
z_{ij}&\mapsto z_{ij} \qquad\text{for other pairs $i$, $j$.}
\end{align*}	

2) Let $g(t)=\exp(t E_{kk})$. We get
\begin{align*}
z_{kj}&\mapsto e^t  z_{kj}, \qquad k<j;\\
z_{ik}& \mapsto e^{-t} z_{ik}, \qquad i<k;\\
z_{ij}&\mapsto z_{ij} \qquad\text{for other pairs $i$, $j$.}
\end{align*}

3) Let $g(t)=\exp (t E_{(k+1)k})$. Then
\begin{align*}
z_{ik}&\mapsto z_{ik}+tz_{i(k+1)}, \qquad i<k;
\\
z_{kj}&\mapsto \frac{z_{kj}}{1+tz_{k(k+1)}},\qquad j>k;
\\
z_{(k+1)m}&\mapsto z_{(k+1)m}+
t\det \begin{pmatrix}
z_{k(k+1)}& z_{km}\\
1&z_{(k+1)m}
\end{pmatrix}, \qquad m>k+1;
\\
z_{ij}&\mapsto z_{ij}\qquad\text{for other pairs $i$, $j$.}
\end{align*}

4) Let $g(t)=\exp(t E_{1n})$, then 
\begin{align*}
 z_{1n}&\mapsto z_{1n}+t,\\
 z_{ij}&\mapsto z_{ij}\qquad\text{for other pairs $i$, $j$.}
\end{align*}

Expressions
(\ref{eq:gener1})--(\ref{eq:gener3}) for the action of the Lie algebra $\frg\frl(n)$
easily follows from these formulas. Also, we get
\begin{equation}
 E_{1n}=\partial_{1n}.
 \label{eq:gener4}
\end{equation}

\sm

{\bf\punct The Jacobian.}

\begin{lemma}
The complex Jacobian $J(Z,g)$ of a transformation 	$Z\mapsto Z^{[g]}$ is
\begin{equation}
J(Z,g)=
\prod_{j=1}^{n-1} \det\bigl([Zg]_{jj}\bigr)^{-2}= 
\prod_{j=1}^{n-1} \det\bigl([Z]_{jn}[Z]_{nj}\bigr)^{-2}.
\label{eq:Jacobian}
\end{equation}
	\end{lemma}

{\sc Proof.} We say that a function 
$$c:N_n^+\times \GL(n,\C)\to \C$$ 
satisfies the {\it chain identity}
if
\begin{equation}
c(Z,gh)=c(Z,g) c(Z^{g},h).
\label{eq:chain}
\end{equation}
The Jacobian satisfies this identity. Let us verify  that the right hand side of
(\ref{eq:Jacobian}) also satisfies it. It suffices to consider one factor 
$\det[Zg]_{jj}^{-2}$. 

Let us evaluate $\det[Z^{[g]}h]_{jj}$.
Represent $Z$, $g$, $h$ etc. as block $j+(n-j)$-matrices,
\begin{equation}
Z=\begin{pmatrix} Z_{11}&Z_{12}\\0&Z_{22}
\end{pmatrix},
\qquad g:=\begin{pmatrix}
G_{11}&G_{12}\\G_{21}&G_{22}
\end{pmatrix},\qquad
h:=\begin{pmatrix}
H_{11}&H_{12}\\H_{21}&G_{22}
\end{pmatrix}.
\label{eq:ZGH}
\end{equation}
Represent $Zg$ as $C^{-1} U$, where $U\in N_{n}^+$, $C\in B_n^-$. Then
$$
\begin{pmatrix}
C_{11}&0\\C_{21}&C_{22}
\end{pmatrix}
\begin{pmatrix}
Z_{11}G_{11}+Z_{12}G_{21}&Z_{11}G_{12}+Z_{12}G_{22}\\
Z_{22}G_{21}&Z_{22}G_{22}
\end{pmatrix}=
\begin{pmatrix} U_{11}&U_{12}\\0&U_{22}
\end{pmatrix}
.
$$
Since $\det U_{11}=1$, we have 
$$\det C_{11}=\det (Z_{11}G_{11}+Z_{12}G_{21})^{-1}.$$
Next,
\begin{multline*}
\det [Z^{[g]}h]_{jj}=\det [Uh]_{jj}
=\det[U_{11}H_{11}+U_{12}H_{21}]
=\\=
\det C_{11} \cdot \det\bigl( (Z_{11}G_{11}+ Z_{12}G_{21})H_{11}+
(Z_{11}G_{12}+ Z_{12}G_{22})H_{21}\bigr)
=\\=
\frac{\det \bigl( Z_{11}(G_{11} H_{11}+G_{12} H_{21})+ Z_{12}(G_{21}H_{21}+G_{})\bigr)}
{\det (Z_{11}G_{11}+Z_{12}G_{21})}
=\frac{\det[Zgh]_{jj}}{\det[Zg]_{jj}}
,\end{multline*}
this proves the desired statement.

Since the both sides of (\ref{eq:Jacobian}) satisfy the chain identity
(\ref{eq:chain}), it suffices to verify the identities for a system
of generators of $\GL(n,\C)$. Thus we can verify (\ref{eq:Jacobian}) for
elements 
one-parametric subgroups $\exp(tE_{k(k+1)})$, $\exp(tE_{kk})$, $\exp(tE_{(k+1)k})$.
Formulas for the corresponding transformations are present in the previous subsection.
Only the case $\exp(tE_{(k+1)k})$ requires a calculation.
In this case,  the Jacobi matrix is triangular. Its diagonal values
are

\sm

--- $(1+t z_{k(k+1)})^{-2}$ for $z_{k(k+1)}$;

\sm

--- $(1+t z_{k(k+1)})^{-1}$ for $z_{km}$, where $m>k+1$;

\sm

--- $(1+t z_{k(k+1)})$ for $z_{(k+1)m}$, where $m>k+1$.

\sm

The Jacobian $(1+t z_{k(k+1)})^{-2}$.

On the other hand, in the product (\ref{eq:Jacobian}) the $j$-th factor
is $(1+t z_{k(k+1)})^{-2}$, other factors are 1.
\hfill $\square$

\sm

{\bf\punct  Proof of Proposition \ref{prop:invariance-of-measure}.}

\sm

{\it The  $\U(n)$-invariance of the $L^2(N_n^+,\mu_\bfp)$-inner product.}
We must check the identity
\begin{equation}
\int_{N_n^+} |(\rho_\bfp(g) \phi)(Z)|^2\,d\mu_\bfp(Z)=
\int_{N_n^+} |\phi(Z)|^2\,d\mu_\bfp(Z).
\label{eq:unitary}
\end{equation}
Fix $g\in\U(n)$.
Substituting $Z=U^{[g]}$ to the right-hand side, we get
$$
\int_{N_n^+} \bigl|\phi(U^{[g]})\bigr|^2 \prod_{j=1}^{n-1} 
\det\left([U^{[g]}]_{jn}\bigl(U^{[g]}]_{jn}\bigr)^* \right)^{-(p_j-p_{j+1})-2}
 \,d\dot{\bigl(U^{[g]}\bigr)}.
$$
In notation (\ref{eq:ZGH}) we have
$$
Ug=\begin{pmatrix}
U_{11}G_{11}+U_{12}G_{21}&U_{11}G_{12}+U_{12}G_{22}\\
*&*
\end{pmatrix}
.
$$
Therefore $[U^{[g]}]_{jn}$ has the form
$$
\begin{pmatrix}
C_{11}(U_{11}G_{11}+U_{12}G_{21})&C_{11}(U_{11}G_{12}+U_{12}G_{22})
\end{pmatrix}
$$
where $C_{11}$ is a lower-triangular matrix, and 
$$
\det\bigl(C_{11}(U_{11}G_{11}+U_{12}G_{21})\bigr)=1.
$$
Hence
\begin{multline*}
\det\left([U^{[g]}]_{jn}\bigl(U^{[g]}]_{jn}\bigr)^*\right)
=\\=
\det\left(\bigl([Ug]_{jn}\bigl([Ug]_{jn}\bigr)^*\right)
\det(U_{11}G_{11}+U_{12}G_{21})^{-1}
\det (U_{11}G_{11}+U_{12}G_{21})^{*-1}.
\end{multline*}
Next, 
\begin{multline*}
\bigl([U]_{jn}\bigl([U]_{jn}\bigr)^*=\\=
(U_{11}G_{11}+U_{12}G_{21})(G_{11}^*U_{11}^*+G_{21}^*U_{12}^*)+
  (U_{11}G_{12}+U_{12}G_{22})(G_{12}^*U_{11}^*+G_{22}^*U_{12}^*)
  =\\=
  U_{11}(G_{11}G_{11}^*+G_{12}G_{12}^*)U_{11}^*+
  U_{11}(G_{11}G_{21}^*+G_{12}G_{22}^*)U_{12}^*+
  \\+
  U_{12}(G_{21}G_{11}^*+G_{22}G_{12}^*)U_{11}^*+
  U_{12}(G_{21}G_{21}+G_{22}G_{22}^*)U_{22}^*
\end{multline*}
The matrix $g$ is unitary, $gg^*=1$. Therefore the last expression equals to
$$
U_{11}\cdot 1 \cdot  U_{11}^*+U_{11}\cdot 0 \cdot U_{12}^*
+U_{12}\cdot 0\cdot U_{11}^*+U_{12}\cdot 1\cdot  U_{12}^*.
$$
Thus
$$
\det\left([U^{[g]}]_{jn}\bigl(U^{[g]}]_{jn}\bigr)^*\right)
=
\det\left(\bigl([U]_{jn}\bigl([U]_{jn}\bigr)^*\right)
\bigr|\det(U_{11}G_{11}+U_{12}G_{21})\bigr|^{-2}.
$$
Keeping in mind
$$
d\dot{\left(U^{[g]}\right)}=\bigr|\det(U_{11}G_{11}+U_{12}G_{21})\bigr|^{-4} \,d\dot Z
,$$
we get
\begin{multline*}
\prod_{j=1}^{n-1} 
\det\left([U^{[g]}]_{jn}\bigl(U^{[g]}]_{jn}\bigr)^* \right)^{-(p_j-p_{j+1})-2}
\,d\dot{\bigl(U^{[g]}\bigr)}
=\\=
\sigma(U,g)
\prod_{j=1}^{n-1} 
\det\left([U]_{jn}\bigl(U]_{jn}\bigr)^*\right)\, d\dot{U}
,
\end{multline*}
where
$$
\sigma(U,g)=
\biggl|\prod_{j=1}^{n-1} \det([Ug]_{jj})^{(p_j-p_{j+1})}\biggr|^2.
$$
Let us represent $Ug$ in the form $D W$, where $D\in B_n^-$, $g\in N_n^+$.
Then
\begin{align}
\det([Ug]_{kk})&=\prod_{i=1}^k b_{ii}(U,g),  \quad k\le n-1,
\label{eq:det-bjj1}
\\
\det(g)&= \prod_{i=1}^n b_{ii}(U,g),
\label{eq:det-bjj2}
\end{align}
here $b_{ii}(U,g)$ are the diagonal elements of $Ug$.
Therefore, 
$$
\sigma(U,g)=\biggl|\prod_{j=1}^{n-1} b_{jj}(U,g)^{p_j-p_{j+1}}\biggr|^2
$$
Since $g\in\U(n)$, we have $|\det(g)|=1$. Keeping 
in  mind
(\ref{eq:det-bjj2}), we come to
$$
\sigma(U,g)=\biggl|\prod_{j=1}^{n} b_{jj}(U,g)^{p_j}\biggr|^2.
$$
Therefore the right-hand side of (\ref{eq:unitary}) equals
$$
\int_{N_n^+} \biggl| f(U^{[g]}) \prod_{j=1}^{n} b_{jj}(U,g)^{p_j}  \biggr|^2
\,d\dot U,
$$
i.e., coincides with the left-hand side.

\sm

Thus the group $\U(n)$ acts in $L^2(N_n^+,d\mu_\bfp)$ by the unitary
operators  (\ref{eq:group-action}). 

\sm

{\it Intersection of the space of $L^2$ and the space of polynomials.}
Denote this intersection by $W$. 

First, $W$ is $\U(n)$-invariant.
Indeed, for $\psi\in W$ we have $\rho_\bfp(g) \psi\in L^2$.
By definition, $\rho_\bfp(g) \psi$ is a rational holomorphic function.
Represent it
as an irreducible fraction $\alpha(Z)/\beta(Z)$.
Let $Z_0$ be a non-singular point of the manifold $\beta(Z)=0$, let $\cO$
be  neighborhood of $Z_0$. It is easy to see that
$$
\int_\cO |\rho_\bfp(g) \psi|^2\,d\dot Z=\infty.
$$
Therefore $\rho_\bfp(g) \psi$ is a polynomial.

\sm

Second, the space $W$ is non-zero. For instance it contains
a constant function. Indeed, $\int_{N_n^+} d\mu_\bfp$
was evaluated in \cite{Ner-hua} and it is finite.

\sm

The third, $W$ is finite-dimensional.
Indeed, our measure has the form $|r(Z)|^{-2}\,d\dot Z$, where
$r(Z,\ov Z)$ is a polynomial in $Z$, $\ov Z$. A polynomial $\phi\in W$
satisfies the condition
$$
\int_{N_n^+}\left|\frac{\phi(Z)}{r(Z,\ov Z)}\right|^2\,d\dot Z<\infty.
$$
Clearly, degree of $\phi$ is uniformly bounded. 

\sm

Next, the operators (\ref{eq:group-action}) determine a unitary
representation in $L^2(N_n^+)$ and this representation is an element
of the principal non-unitary series (see, e.g., \cite{GN} or  \cite{Zhe-book}, Addendum).
But a representation of the principal series can not have more than one finite-dimensional
subrepresentation.
\hfill $\square$

\sm

\section{The formula for the kernel}

\COUNTERS

In this section we prove  Proposition \ref{prop:kernel}.

\sm

{\bf \punct Functions $\Phi[\dots]$, $\Psi[\dots]$.}
Consider the matrices $Z$, $Z^{cut}$, $U$, $U^{ext}$ as in Subs. \ref{ss:normalization}.
We use notation (\ref{eq:submatrix-1}) for their sub-matrices.
Let $\alpha=1,\dots,n-1$.
Let $I$, $J$ be sets of integers,
$$I: \,\,0<i_1<\dots <i_{n+1-\alpha}\le n+1,
\qquad
J:\,\,0<j_1<\dots<j_\alpha\le n
.
$$
Denote
\begin{equation}
\Phi_\alpha[I;J]=\Phi_\alpha[i_1,\dots,i_{n+1-\alpha};j_1,\dots j_\alpha]:=
\det \begin{pmatrix}
[Z(I)]\\
[U^{ext}(J)]
\end{pmatrix}
.
\end{equation}
Let $I$, $J$ be sets of integers,
$$I: \,\,0<i_1<\dots <i_{n-\alpha}\le n+1,
\qquad
J:\,\,0<j_1<\dots<j_\alpha\le n.
$$
Denote
\begin{equation}
\Psi_\alpha[I;J]=\Psi_\alpha[i_1,\dots,i_{n-\alpha};j_1,\dots j_\alpha]:=
\det \begin{pmatrix}
[Z^{cut}(I)]\\
[U(J)]
\end{pmatrix}
.
\end{equation}

We have two collections of variables $z_{ij}$, $u_{i'j'}$.
Zhelobenko operators $R_{kl}$  in $z$ and $u$ we denote 
by
$$
R^z_{kl},\qquad R^u_{kl}.
$$

\begin{lemma}
	\label{l:R-Phi}
\begin{equation}
R_{kl}^z\Phi_\alpha[I,J]=\begin{cases}
0, &\text {if $k\notin I$};
\\
0, &\text {if $k\in I$, $l\in I$};
\\
(-1)^{\theta(I,k,l)} \Phi_{\alpha}[I^\circ;J],&\text{if  $k\in I$, $l\notin I$}
.
\end{cases}
\end{equation}
where $I^\circ$ is obtained from $I$ by replacing of $k$ by $l$, and the corresponding
change of order. If $k=i_s$ and $i_t<l<i_{t+1}$, then
$$
I^\circ=(i_1,\dots,i_{s-1}, i_{s+1},\dots, i_{t-1}, l, i_{t+1}, \dots, i_{n+1-\alpha})
$$
and 
$$\theta(I,k,l)=t-s.$$

The same properties {\rm(}with obvious modifications{\rm)} hold for
$$
R^u_{kl} \Phi_\alpha[I;J],\qquad R^z_{kl} \Psi_\alpha[I;J],\qquad R^u_{kl} \Psi_\alpha[I;J].
$$
\end{lemma}

{\sc Proof.} Recall that 
$$R^z_{kl}=\partial_{kl}+\sum_{j>l} z_{lj}\partial_{kj}.$$
If $k\notin I$, then $\Phi_\alpha[I,J]$ does not depend on variables
$z_{km}$. Therefore we get 0. If $k\in I$, these variables are present only in
one row. Application of  $R^z_{kl}$ is equivalent
to a change of the row
$$\begin{pmatrix}
0&\dots&0&1&z_{k(k+1)}&\dots&z_{kl}&z_{k(l+1)}&\dots&z_{k(n+1)}
\end{pmatrix}
$$
of the matrix 
$\begin{pmatrix}
[Z(I)]\\
[U^{ext}(J)]\end{pmatrix}$ by the row 
$$\begin{pmatrix}
0&\dots&0&0&0&\dots&1&z_{l(l+1)}&\dots&z_{k(n+1)}
\end{pmatrix}.
$$
If $l\in I$, then  this row is present in the initial matrix, and  again
$\det=0$. Otherwise, we come to $\Phi_\alpha[I^\circ;J]$
up to the order of rows.
\hfill $\square$ 

\sm

\begin{corollary}
	\label{cor:R-Phi-1}
For any $k$, $l$, $\alpha$,
\begin{align*}
(R^z_{kl})^2\Phi_\alpha[I;J]=0, \quad (R^u_{kl})^2\Phi_\alpha[I;J]=0,
\\
(R^z_{kl})^2\Psi_\alpha[I;J]=0, \quad (R^u_{kl})^2\Psi_\alpha[I;J]=0.
\end{align*}
\end{corollary}

\begin{corollary}
	\label{cor:R-Phi-2}
	\begin{equation}
R^z_{kl}	R^z_{lm} \Phi_\alpha[I;J]=
\begin{cases}
- R^z_{km}\Phi_\alpha[I;J],&\quad\text{if $k\in I$, $l\notin I$, $m\notin I$;}
\\
0,&\quad \text{otherwise}. 
\end{cases}
	\end{equation}
The same property holds for operators $R^u$ and for $\Psi_\alpha[I;J]$.
\end{corollary}

{\bf\punct Verification of the Zhelobenko conditions.}

\begin{lemma}
$L_{\bfp,\bfq}\in V_\bfp\otimes V_{\bfq}$.	
\end{lemma}

{\sc Proof.} By the interlacing conditions, 
$L_{\bfp,\bfq}$ is a polynomial.
 We must verify the identities
\begin{equation}
(R_{k(k+1)}^z)^{p_k-p_{k+1}+1} L=0,\qquad 
(R_{k(k+1)}^u)^{q_k-q_{k+1}+1} L=0.
\label{eq:L-in-tensor}
\end{equation}
To be definite verify the first equality,
\begin{multline*}
(R^z_{k(k+1)})^{p_k-p_{k+1}+1} L_{\bfp,\bfq}=\\= 
(R^z_{k(k+1)})^{p_k-p_{k+1}+1}
 \left(\Psi_{n-k}^{-p_{k+1}-q_{n+1-k}}\Phi_{n+1-k}^{p_k+q_{n+1-k}} 
\cdot 
\biggl\{\text{remaining factors} \biggl\}
\right)
=\\= 
(R^z_{k(k+1)})^{p_k-p_{k+1}+1}
\left(\Psi_{n-k}^{-p_{k+1}-q_{n+1-k}}\Phi_{n+1-k}^{p_k+q_{n+1-k}}\right) 
\cdot 
\biggl\{\text{remaining factors} \biggl\}
\end{multline*}
The sum of exponents (the both exponents are positive) is
$$
\bigl(-p_{k+1}-q_{n+1-k}\bigr)+\bigl(p_k+q_{n+1-k}\bigr)=p_k-p_{k+1},
$$
and we obtain the desired 0. If $k=n$ the factor $\Psi_{n-k}$ is absent,
we have
$$
(R^z_{n(n+1)})^{p_n-p_{n+1}+1}\left(\Phi_1^{p_n+q_1} \right)
\cdot 
\biggl\{\text{remaining factors} \biggl\}.
$$
By the interlacing condition $-p_{n+1}\ge q_1$ and we
again obtain 0. 
\hfill $\square$

\sm

{\bf\punct  Invariance of the kernel.} 
Consider matrices $\Phi_\alpha$, $\Psi_\alpha$  defined by (\ref{eq:Phi-Psi}).

\begin{lemma}
	Let $g\in\GL(n,\C)$. Denote $\wt g=\begin{pmatrix} g&0\\0&1\end{pmatrix}\in \GL(n+1,\C)$.
	Then
\begin{multline}
\Phi(Z^{[\wt g]},U^{[g]})=
\\=
\Phi(Z,U) 
\det(g)\det\bigl([Z\wt g]_{(n+1-\alpha)(n+1-\alpha)}\bigr)^{-1}
\det \bigl([Ug]_{\alpha\alpha}\bigr)^{-1}.
\label{eq:Phi-g}
\end{multline}
\begin{multline}
\Psi(Z^{[\wt g]},U^{[g]})=\\=
\Psi(Z,U) 
\det(g)\det\bigl([Z\wt g]_{(n-\alpha)(n-\alpha)}\bigr)^{-1}
\det \bigl([Ug]_{\alpha\alpha}\bigr)^{-1}.
\label{eq:Psi-g}
\end{multline}
\end{lemma}

{\sc Proof.}  We write the Gauss decompositions of $Z\wt g$, $Ug$,
$$
Z\wt g=B P,\qquad Ug= C Q
.
$$
Then the Gauss decomposition of $U^{ext}$ (see (\ref{U-ext})) is
$$
U^{ext}\wt g= \begin{pmatrix} C&0\\0&1\end{pmatrix}\begin{pmatrix} Q&0\\0&1\end{pmatrix}
,$$ 
denote factors in the right-hand side by $\wt C$, $\wt Q$. Then
\begin{multline*}
\begin{pmatrix}
[Z^{[\wt g]}]_{(n+1-\alpha)(n+1)}
\\
[(U^{ext})^{[\wt g]}]_{\alpha(n+1)}
\end{pmatrix}
=\begin{pmatrix}[P]_{(n+1-\alpha)(n+1)}\\
[\wt Q]_{\alpha(n+1)} \end{pmatrix}
=\\=
\begin{pmatrix}
([\wt B]_{(n+1-\alpha)(n+1-\alpha)})^{-1}&0\\
0& ([\wt C]_{\alpha\alpha})^{-1}
\end{pmatrix}
\begin{pmatrix}
[Z]_{(n+1-\alpha)(n+1)}
\\
 [\wt U]_{\alpha(n+1)}
\end{pmatrix}
\cdot g
\end{multline*}
We pass to determinants in the left-hand side and the right-hand side.
Keeping in  mind
\begin{align*}
&\det [\wt B]_{(n+1-\alpha)(n+1-\alpha)}= \det[\wt Z\wt g]_{(n+1-\alpha)(n+1-\alpha)},
\\
&\det [\wt C]_{\alpha\alpha}= \det[U^{ext} \wt g]_{\alpha\alpha}=
\det[U g]_{\alpha\alpha},
\end{align*}
we come to (\ref{eq:Phi-g}). Proof of (\ref{eq:Psi-g}) is similar.
\hfill $\square$

\sm

{\sc Proof of Proposition \ref{prop:kernel}.}
Applying the lemma and (\ref{eq:det-bjj1}), we get
$$L_{\bfp,\bfq}(Z^{[\wt g]},U^{[g]})
=L_{\bfp,\bfq}(Z,U) \prod_{j=1}^n b_{jj}(Z,g)^{-q_1-p_j}
\prod b_{jj}(U,g)^{-p_n-q_j}
\det(g)^{p_1+q_n} .
$$
Since
$$
b_{(n+1)(n+1)}(Z,\wt g)=1, \qquad
\det(g)=\prod_{j=1}^{n} b_{jj}(Z,g)=\prod_{j=1}^{n} b_{jj}(U,g)
,
$$
we come to
$$
L_{\bfp,\bfq}(Z^{[\wt g]},U^{[g]})
=L_{\bfp,\bfq}(Z,U) \prod_{j=1}^n b_{jj}(Z,g)^{-p_j}
\prod b_{jj}(U,g)^{-q_j}.
$$
Therefore,
$$
(\rho_\bfp(\wt g)\otimes \rho_\bfq(g)) L_{\bfp,\bfq}=L_{\bfp,\bfq}.
$$
Thus $ L_{\bfp,\bfq}$ is a non-zero $\GL(n,\C)$-invariant vector
in $V_\bfp\otimes V_\bfq$. We also know that such a vector is unique up to a 
scalar factor.
\hfill $\square$

\section{The calculation}

\COUNTERS 

Below $\bfq\in \Lambda_n$, $\bfr=\bfp^*\in \Lambda_{n+1}$ are signatures.
The signatures $\bfr$ and $\bfq$ are interlacing.
The kernel $L$ is the same as above (\ref{eq:kernel}),
$$
L:=L_{\bfq}^\bfr (Z,U)=	\Phi_1^{q_1-r_2}\Psi_1^{r_2-q_2} \Phi_2^{q_2-r_3} 
\Psi_2^{r_3-q_3}
\dots \Phi_n^{q_n-r_{n+1}}.
$$
Below 
$$R_{kl}:=R_{kl}^u.$$

We must verify the differential-difference equations
$$
E_{1(n+1)} L= F_{1(n+1)}L,\qquad E_{(n+1)n} L= F_{(n+1)n} L
,$$
 $F_{1(n+1)}$, $F_{(n+1)n}$ are  given by (\ref{eq:formula-1})--(\ref{eq:formula-1}).
Evaluations of the left-hand sides is easy, it is contained in the next subsection.
Evaluation of the right-hand sides is  tricky and occupies the remaining part of the section.

\sm

{\bf \punct Evaluation of $E_{1(n+1)} L$.}

\begin{lemma}
	\begin{equation}
	 E_{1(n+1)} L=
\Bigl(-\frac {q_1-r_2}{\Phi_1}+ 	\sum_{\alpha=2}^{n}
\frac
{(q_\alpha-r_{\alpha+1})R_{1\alpha} \Psi_{\alpha-1}}{\Phi_\alpha}
	\Bigr) \cdot L.
	\label{eq:E1(n+1)L}
	\end{equation}
\end{lemma}

 {\sc Remark.}
For $\alpha=1$ we immediately get 
$\frac{\partial \Phi_1}{\partial z_{1(n+1)}}=-1$.
 On the other hand,
taking formally $\alpha=1$ in the sum in (\ref{eq:E1(n+1)L}) we get
$R_{11} \Psi_0/ \Phi_1$.
At first glance,  we must assume that $\Psi_0=1$ and that $R_{11}$ is the  identical operator.
Under this assumption, $R_{11} \Psi_0/ \Phi_1=1/\Phi_1$.  However,
this gives an incorrect sign.
\hfill $\boxtimes$

\sm

In the following two proofs we need manipulations with determinants.
To avoid huge matrices or compact notation, which are difficult for reading,
we expose calculations for matrices having a minimal size that allows to
visualize  picture.

\sm

{\sc Proof.} We have (see (\ref{eq:gener4}))
$$
E_{1(n+1)} L=\frac{\partial L}{\partial z_{1(n+1)}}
=
\Bigl(\sum_{\alpha=1}^{n} (q_\alpha-r_{\alpha+1})
 \frac{\partial \Phi_\alpha}{\partial z_{1(n+1)}}\cdot\Phi_\alpha^{-1} \Bigr)
	\cdot L.
	$$
 Take $n=3$ and $\alpha=2$,
 \begin{equation}
 \Phi_2=
 \det\begin{pmatrix}
 	1&z_{12}& z_{13}& z_{14}\\
 	0& 1& z_{23}& z_{24}\\
 	1&u_{12}&u_{13}&0\\
 	0&1&u_{23}&0
 \end{pmatrix}
 .
 \label{eq:Phi2}
 \end{equation}
 Then
 $$
 \frac{\partial \Phi_2}{\partial z_{14}}
 = - \det\begin{pmatrix}
 0& 1& z_{23}\\
 1&u_{12}&u_{13}\\
 0&1&u_{23}
 \end{pmatrix}
 =  \det\begin{pmatrix}
 1& z_{23}\\
 1&u_{23}
 \end{pmatrix}=
 \det\begin{pmatrix}
 	1&z_{12}& z_{13}\\
 0& 1& z_{23}\\
 0&1&u_{23}
 \end{pmatrix}=R_{12}\Psi_1.
 $$

\begin{lemma}
$$
E_{(n+1)n} L=
\Bigl(z_{n(n+1)}(r_1-q_1)-\sum_{j=1}^{n-1}
	\frac{(r_{\alpha+1}-q_{\alpha+1})R_{(\alpha+1)n}\Phi_{\alpha+1}}{\Psi_\alpha}\Bigr)
	\cdot L.
$$	
\end{lemma}

{\sc Proof.} Split $E_{(n+1)n}$ as a sum 
$$
E_{(n+1)n}=
D+(r_1-r_2)z_{n(n+1)},\qquad
D:=\sum_{j=1}^{n-1}z_{j(n+1)}\partial_{jn}-z_{n(n+1)}^2 \partial_{n(n+1)}
.
$$
Obviously,
\begin{multline*}
E_{(n+1)n} L=
\\=
\Bigl((r_1-r_2)z_{n(n+1)}+\sum_{\alpha=1}^{n}
\frac{(q_\alpha-r_{\alpha+1})D\Phi_\alpha}{\Phi_\alpha}+
\sum_{\beta=1}^{n-1} \frac{(r_{\beta+1}-q_{\beta+1})D\Psi_\beta}{\Psi_\beta}\Bigr)\cdot L
\end{multline*}
Let us evaluate all $D\Phi_\alpha$, $D\Psi_\alpha$.

\sm

1) $D\Phi_\alpha=0$ for all $\alpha\ne 1$.
Take $n=3$, $\alpha=2$,
$$
D  \Phi_2= D
 \det\begin{pmatrix}
 	1&z_{12}& z_{13}& z_{14}\\
 	0& 1& z_{23}& z_{24}\\
 	1&u_{12}&u_{13}&0\\
 	0&1&u_{23}&0
 \end{pmatrix}
=
\det\begin{pmatrix}
1&z_{12}& z_{14}& z_{14}\\
0& 1& z_{24}& z_{24}\\
1&u_{12}&0&0\\
0&1&0&0
\end{pmatrix}
=0
.$$
Variables $z_{jn}$ (in our case $z_{j4}$) are present only in the
next-to-last column. We apply the operator $D$ to this column and come to
a matrix with coinciding columns. 

\sm

2) $D\Phi_1=-z_{n(n+1)}\Phi_1$. 
Take $n=3$,
\begin{multline*}
	D\Phi_1=D
	\det
	\begin{pmatrix}
		1&z_{12}&z_{13}&z_{14}\\
		0&1     &z_{23}&z_{24} \\
		0&0     &1     &z_{34}\\
		1&u_{12}&u_{13}&0
	\end{pmatrix}=\\=
\det
		\begin{pmatrix}
			1&z_{12}&z_{14}&z_{14}\\
			0&1     &z_{24}&z_{24} \\
			0&0     &0     &z_{34}\\
			1&u_{12}&0     &0
		\end{pmatrix}
+	z_{34}^2	\det
		\begin{pmatrix}
			1&z_{12}&z_{13}\\
			0&1     &z_{23} \\
		    1&u_{12}&u_{13}
		\end{pmatrix}
		=\\=
-z_{34}\det\begin{pmatrix}
1&z_{12}&z_{14}\\
0&1     &z_{24} \\
1&u_{12}&0
\end{pmatrix}
+	z_{34}^2	\det
\begin{pmatrix}
1&z_{12}&z_{13}\\
0&1     &z_{23} \\
1&u_{12}&u_{13}
\end{pmatrix}		
=\\=
-z_{34}\det\begin{pmatrix}
1&z_{12}&z_{13}&z_{14}\\
0&1     &z_{23}&z_{24} \\
0&0     &1     &0  \\
1&u_{12}&u_{13}&0
\end{pmatrix}
-	z_{34}	\det
\begin{pmatrix}
1&z_{12}&z_{13}&0\\
0&1     &z_{23}&0 \\
0&0     &1    &z_{34}  \\
1&u_{12}&u_{13}&0
\end{pmatrix}
=\\=
-z_{34}\det\begin{pmatrix}
1&z_{12}&z_{13}&z_{14}\\
0&1     &z_{23}&z_{24} \\
0&0     &1     &z_{34}  \\
1&u_{12}&u_{13}&0
\end{pmatrix}
=-z_{34}\Phi_1
.
\end{multline*}

3) $D\Psi_\alpha=-R_{(\alpha+1)n}\Phi_{\alpha+1}$.
Take $n=4$, $\alpha=1$.
\begin{multline*}
D \Psi_1=
D\det\begin{pmatrix}
1&z_{12}&z_{13}&z_{14}\\
0&1     &z_{24}&z_{24}\\
0&0     &1     &z_{34}\\
1&u_{12}&u_{13} &u_{14}
\end{pmatrix}=
\det\begin{pmatrix}
1&z_{12}&z_{13}&z_{15}\\
0&1     &z_{24}&z_{25}\\
0&0     &1     &z_{35}\\
1&u_{12}&u_{13} &0
\end{pmatrix}
=\\=-\det
\begin{pmatrix}
1&z_{12}&z_{13}&z_{14}&z_{15}\\
0&1     &z_{24}&z_{24}&z_{25}\\
0&0     &1     &z_{25}&z_{35}\\
1&u_{12}&u_{13}&u_{14}&0    \\
0&0     &0     &1     &0
\end{pmatrix}=-R_{24}\Phi_2.
\end{multline*}

Thus we know all $D\Phi_\alpha$, $D\Psi_\alpha$. This gives us the desired statement.
\hfill $\square$

\sm

{\bf \punct Quadratic relations.}
Now we wish to write a family of quadratic relations between different functions of type
$R_{kl}\Phi_\alpha$, $R_{kl}\Psi_\beta$. We introduce notation
\begin{equation}
 R_{kk}\Phi_{k-1}:=+\Phi_{k-1},\qquad \boxed{R_{kk}\Psi_{k-1}=-\Psi_{k-1},\quad
 R_{11}\Psi_0=-\Psi_0=1}\,
.\end{equation}
We do not introduce operators $R_{kk}$ and $ R_{kk}\Phi_{k-1}$,  $R_{kk}\Psi_{k-1}$
are only symbols used in formulas. In particular,  formula
(\ref{eq:E1(n+1)L}) now can be written the form
	\begin{equation*}
	 E_{1(n+1)} L=
L\cdot \sum_{\alpha=1}^{n}
\frac
{(q_\alpha-r_{\alpha+1})R_{1\alpha} \Psi_{\alpha-1}}{\Phi_\alpha}
		\end{equation*}
		without a term with abnormal sign. 
Below it allows to avoid numerous anomalies in the formulas, duplications of formulas
and branchings of calculations.

\begin{lemma}
{\rm a)}
Let $m<\alpha\le \beta$, 
\begin{equation}
 \sum_{j=\alpha+1}^{\beta+1} R_{mj} \Phi_\alpha \cdot R_{j(\beta+1)}\Psi_\beta= -\Phi_\alpha \cdot R_{m(\beta+1)}
  \Psi_\beta
  +
  \Phi_{\beta+1}\cdot R_{m\alpha}\Psi_{\alpha-1}.
  \label{eq:Plu-1}
  \end{equation}

	{\rm b)} Let $m<\alpha<\beta$. Then
	\begin{equation}
	\sum_{j=\alpha+1}^{\beta+1}
	 R_{mj}\Psi_\alpha\cdot R_{j(\beta+1)}\Psi_\beta=-\Psi_\alpha\cdot R_{m(\beta+1)}\Psi_\beta
	 .
	  \label{eq:Plu-2}
	 \end{equation}
	

	{\rm c)} Let $m<\alpha<\beta$. Then 
\begin{equation}
 \sum_{j=\alpha+1}^\beta R_{m j}\Phi_\alpha\cdot R_{j(\beta+1)} \Phi_\beta
 =-\Phi_\alpha \cdot R_{m(\beta+1)}-R_{m(\beta+1)}\Phi_\alpha\cdot \Phi_\beta
 .
  \label{eq:Plu-3}
\end{equation}

\end{lemma}

{\sc Proof.} Denote by $\xi$ the row $\begin{pmatrix}0&0&\dots&1 \end{pmatrix}$
of the length $(n+1)$.
We write a matrix $\Delta:=\begin{pmatrix}Z \\U^{ext}\\\xi \end{pmatrix}$.
Let us enumerate $z$-rows of these matrix by 1, 2, 3, \dots  and $u$-rows
by marks $\ov 1$, $\ov 2$, $\ov 3$, \dots. 
Consider the following minors of $\Delta$:
\begin{align}
\Phi_\alpha&=\det\Delta[1,\dots,n+1-\alpha;\ov 1,\dots,\ov \alpha]
;
\nonumber
\\
R_{m(\beta+1)} \Psi_\beta&
=(-1)^{\beta-m}
\det \Delta[1,\dots,n-\beta;\ov 1,\dots, \widehat{\ov m},\dots,  \beta+1;\xi]
;
\label{eq:pplu-1}
\end{align}
Notice that we  $\xi$ in the last row allows to replace the last
determinant by 
$$
\det\Bigl[ \Delta[1,\dots,n-\beta;\ov 1,\dots, \widehat{\ov m},\dots,  \beta+1;\xi]\Bigr]_{n n}
.
$$
and this proves (\ref{eq:pplu-1}).

\sm

We wish to apply one of the Pl\"ucker identities to the product of these minors
(see, e.g., \cite{Ful}, Section 9.1). We take $\ov m$ from the collection
$(1,\dots,n+1-\alpha;\ov 1,\dots,\ov \alpha)$, exchange it with an element of the collection
$(1,\dots,n-\beta;\ov 1,\dots, \widehat{\ov m},\dots,  \beta+1;\xi)$
and consider the product of the corresponding minors.
Next, we take the sum of all such products of $\Delta$. In this way, we obtain
{\small
\begin{align}
&\det \Delta[1,\dots,n+1-\alpha;\ov 1,\dots,\ov \alpha]\cdot
\det \Delta[1,\dots,n-\beta;\ov 1,\dots, \widehat{\ov m},\dots,  \ov{\beta+1};\xi]
\label{eq:phipsi-0}
=\\=
&\det \Delta[1,\dots,n+1-\alpha;\ov 1,\dots,\ov{m-1}, 1,\ov{m+1}, \dots,  \ov \alpha]\times
\label{eq:phipsi-1}
\\&\qquad\qquad\qquad\qquad
\times
\det \Delta[\ov m, 2,\dots,n-\beta;\ov 1,\dots, \widehat{\ov m},\dots,  \ov{\beta+1};\xi]+\dots+
\label{eq:phipsi-2}
\\
+&\det \Delta[1,\dots,n+1-\alpha;\ov 1,\dots,\ov{m-1}, \ov 1,\ov{m+1}, \dots,  \ov \alpha]\times
\label{eq:phipsi-3}
\\
&\qquad\qquad\qquad\qquad
\times
\det \Delta[1,\dots,n-\beta;\ov m,\ov2,\dots, \widehat{\ov m},\dots,  \ov{\beta+1};\xi]+\dots+
\label{eq:phipsi-4}
\\
&+\det \Delta[1,\dots,n+1-\alpha;\ov 1,\dots,\ov{m-1}, \ov  {\alpha+1},\ov{m+1}, \dots,  \ov \alpha]\times
\label{eq:phipsi-5}
\\&\qquad
\times
\det \Delta[1,\dots,n-\beta;\ov 1,\dots, \widehat{\ov m},\dots,\ov{\alpha},\ov m,\ov  {\alpha+2},\dots,   \beta+1;\xi]
+\dots+
\label{eq:phipsi-6}
\\
& +\det \Delta[1,\dots,n+1-\alpha;\ov 1,\dots,\ov{m-1}, \ov  {\beta},\ov{m+1}, \dots,  \ov \alpha]\times
\label{eq:phipsi-7}
\\&\qquad\qquad\qquad\qquad
\times
\det \Delta[1,\dots,n-\beta;\ov 1,\dots, \widehat{\ov m},\dots,\ov{\beta-1},\ov m,  \beta+1;\xi]
\label{eq:phipsi-8}
\\
&+\det \Delta[1,\dots,n+1-\alpha;\ov 1,\dots,\ov{m-1}, \ov  {\beta+1},\ov{m+1}, \dots,  \ov \alpha]\times
\label{eq:phipsi-9}
\\
&\qquad\qquad\qquad\qquad\qquad\qquad\det \Delta[1,\dots,n-\beta;\ov 1,\dots, \widehat{\ov m},\dots,  \ov \beta,\ov m;\xi]
+
\label{eq:phipsi-10}
\\
&+\det \Delta[1,\dots,n+1-\alpha;\ov 1,\dots,\ov{m-1}, \xi,\ov{m+1}, \dots,  \ov \alpha]\times
\label{eq:phipsi-11}\\
&\qquad\qquad\qquad\qquad
\det \Delta[1,\dots,n-\beta;\ov 1,\dots, \widehat{\ov m},\dots,  \ov{\beta+1};\ov m]
.
\label{eq:phipsi-12}
\end{align}
}
Look to the summands of our expression up to signs.

\sm

1) After exchanging of $\ov m$ with $1$, \dots, $n-\beta$ we get a matrix $\Delta[\dots]$
with coinciding rows and therefore we obtain 0, see (\ref{eq:phipsi-1}).

\sm

2) After exchanging of $\ov m$ with $\ov 1$, \dots, $\alpha$ we again obtain 0,
see (\ref{eq:phipsi-3}).

\sm

3) The sum
(\ref{eq:phipsi-5})--(\ref{eq:phipsi-8}) corresponds to the sum in the left hand side of the the desired
identity (\ref{eq:Plu-1}) with $j<\beta+1$. 

\sm

4)
The summand (\ref{eq:phipsi-9})--(\ref{eq:phipsi-10})
corresponds to the last term of the sum.

\sm

5) The summand (\ref{eq:phipsi-11})--(\ref{eq:phipsi-12}) corresponds to the expression 
$\Phi_{\beta+1}\cdot R_{m\alpha} \Psi_{\alpha-1}$.

\sm

6) The term (\ref {eq:phipsi-0}) is $\Phi_\alpha\cdot R_{m(\beta+1)} \Psi_{\beta}$.

\sm

 Next,
let us watch signs. Denote by $(i_1i_2\dots i_k)$ a cycle in a substitution.
Denote by $\sigma(\cdot)$ the parity of a substitution. 

Examine summands (\ref{eq:phipsi-5})--(\ref{eq:phipsi-8}).
For $j=\alpha+1$, \dots, $\beta+1$,
we have 
\begin{multline*}
\det \Delta[1,\dots,n+1-\alpha;\ov 1,\dots,\ov{m-1}, \ov  {j},\ov{m+1}, \dots,  \ov \alpha]
=\\=
\sigma\bigl(j\, (m+1)\, \dots\, \alpha \bigr)R_{mj} \Phi_\alpha
= (-1)^{\alpha-m} R_{mj} \Phi_\alpha
\end{multline*}
For $j=\alpha+1$, \dots, $\beta$ we have
\begin{multline*}
\det \Delta[1,\dots,n-\beta;\ov 1,\dots, \widehat{\ov m},\dots,\ov{j-1},\ov m,\ov  {j+1},\dots,   \beta+1;\xi]
=\\=\sigma\bigl(m\dots (j-1)\bigr)\cdot \sigma\bigl((j+1)\dots(\beta+1) \bigr) R_{j(\beta+1)} \Psi_\beta]
=(-1)^{\beta-m+1}
R_{j(\beta+1)} \Psi_\beta].
\end{multline*}
Next, in (\ref{eq:phipsi-9}) we have
$$\det \Delta[1,\dots,n-\beta;\ov 1,\dots, \widehat{\ov m},\dots,  \ov \beta,\ov m;\xi]
=(-1)^{\beta-m} \Psi_\beta
$$
and this gives the anomaly of a sign mentioned above.

Examine the last summand (\ref{eq:phipsi-11})--(\ref{eq:phipsi-12}). We get
\begin{multline*}
 \det \Delta[1,\dots,n+1-\alpha;\ov 1,\dots,\ov{m-1}, \xi,\ov{m+1}, \dots,\ov{\alpha-1},  \ov \alpha]
 =\\
=- \det \Delta[1,\dots,n+1-\alpha;\ov 1,\dots,\ov{m-1}, \ov \alpha,\ov{m+1}, \dots,\ov{\alpha-1},   \xi]=
\\
-\det \Bigl[\Delta[1,\dots,n+1-\alpha;\ov 1,\dots,\ov{m-1}, \ov \alpha,\ov{m+1}, \dots,\ov{\alpha-1}]\Bigr]_{nn}
=-R_{m\alpha}\Psi_{\alpha-1}
\end{multline*}
and
\begin{multline*}
\det \Delta[1,\dots,n-\beta;\ov 1,\dots, \widehat{\ov m},\dots,  \ov{\beta+1};\ov m]
=\\=\sigma\bigl(m\,\dots\, (\beta+1)\bigr) \Phi_{\beta+1}=(-1)^{\beta+1-m}
\Phi_{\beta+1}.
\end{multline*}
Combining all signs we get the desired identity.

\sm

{\sc Proof of the statement b.}
We compose a new matrix $\Delta=\begin{pmatrix}
                                 Z^{cut}\\U
                                \end{pmatrix}$,
   consider the following minors of $\Delta$,
  \begin{align*}
   \Psi_\alpha&=\det \Delta[1,\dots,n-\alpha;\ov 1,\dots, \ov \alpha],
   \\
   R_{m(\beta+1)}\Psi_\beta&=(-1)^{\beta-m}\det\Delta[1,\dots,n-\beta;\ov1, \dots,\wh{\ov m},\dots,\ov\beta,\ov{\beta+1}],
  \end{align*}
and write the Pl\"ucker identity exchanging the  $\ov m$-th row of $\Psi_\alpha$ with 
rows of $R_{m(\beta+1)}\Psi_\beta$. 

\sm

{\sc Proof of the statement c.} We compose a matrix
$\Delta:=\begin{pmatrix}
          Z\\U^{ext}
         \end{pmatrix}$, take the minors
           \begin{align*}
   \Phi_\alpha& =\det \Delta[1,\dots,n+1-\alpha;\ov 1,\dots, \ov \alpha],
   \\
   R_{m(\beta+1)}\Phi_\beta&=(-1)^{\beta-m}\det\Delta[1,\dots,n+1-\beta;\ov1, \dots,\wh{\ov m},\dots,\ov\beta,\ov{\beta+1}].
  \end{align*}
  and write the Pl\"ucker identities exchanging $\ov m$-th of $\Phi_\alpha$ with
  rows of the second minor.
  \hfill $\square$

\sm

\sm

{\bf \punct A recurrence formula.%
\label{ss:requerence}}
Our next purpose is to derive a formula for $F_{1(n+1)}L$.
Let us evaluate one of summands of this expression,
	\begin{equation}
	\Theta_{1n}:=\Bigl(\sum_{I\triangleleft [1,n]}
	\prod_{l\in [1,n]\setminus I} (l-n+q_l-q_n) \cdot R_I
	\Bigr) T^-_n L.
	\label{eq:summand}
	\end{equation}

Recall that the notation $J\triangleleft\{1,2,\dots n\}$ means
that $J$ is a subset in $\{1,2,\dots n\}$ containing $1$ and $n$,
see Subsection \ref{ss:ao}.

\sm

{\sc Remark.} By definition, the right-hand side is a rational expression of the form
$$
\frac{F}{\prod \Phi_\alpha^{s_\alpha}\prod\Psi_\beta^{t_\beta}}.
$$
It is easy to observe (see below) that all the exponents $s_\alpha$, $t_\beta$ are equal 1.
In a calculation below we decompose this function as a sum of 'prime fractions' with denominators
$\Phi_\alpha$, $\Psi_\beta$ and finally get an unexpectedly simple expression.
\hfill $\boxtimes$

\sm

Let $1\le k<l\le n$, $\mu<n$. Denote
\begin{align}
&W^\mu_{kl}:=\frac{(q_{\mu}-r_{\mu+1})R_{kl}\Phi_\mu}{\Phi_\mu}+
\frac{(r_{\mu+1}-q_{\mu+1})R_{kl}\Psi_\mu}{\Psi_\mu}
;
\label{eq:W}
\\
&w^\mu_{kl}:=(-1)^{l-k-1}\prod_{j=k+1}^\mu (q_n-q_j+j-n-1)  \prod_{j=\mu+1}^{l-1} (q_n-q_j+j-n)
\label{eq:w};
\\
&\zeta_{kl}=\sum_{\mu=k}^{l-1} w^\mu_{kl} W^\mu_{kl}+ w_{k(n-1)}^{n-1}\frac{R_{kl}\Psi_{n-1}}{\Psi_{n-1}}.
\label{eq:zeta}
\end{align}	
Notice that the additional  term in the last row is 0 if $l\ne n$.

\begin{lemma}
\begin{align}
&	R_{kl} L= L\cdot  \sum_{\mu: k\le\mu <l} W_{kl}^\mu;
		\label{eq:RL-1}
	\\
&	R_{kl} T_n^- L= L\cdot  \Bigl(\sum_{\mu: k\le\mu <n} W_{kl}^\mu\Bigr)\cdot \frac{\Psi_{n-1}}{\Phi_n},
	\qquad \text{for $l<n$};
	\label{eq:RL-2}
	\\	
&		R_{kn} T_n^-  L= L\cdot  \Bigl(\sum_{\mu: k\le\mu <l} W_{kl}^\mu +
		\frac{R_{kn}\Psi_{n-1}}{\Psi_{n-1}}\Bigr)\cdot \frac{\Psi_{n-1}}{\Phi_n}.
		\label{eq:RL-3}
\end{align}	
\end{lemma}	

{\sc Proof.}   Obviously,
$$
R_{kl}L=L\Bigl(\sum_{j=1}^n \frac{(q_{\mu}-r_{\mu+1})R_{kl  \Phi_\mu}}{\Phi_\mu} +
\sum_{j=1}^n \frac{(r_{\mu+1}-q_{\mu+1})R_{kl  \Psi_\mu}}{\Psi_\mu} \Bigr).
$$
 By Lemma \ref{l:R-Phi}, only terms with $k\le\mu <l$ give a nonzero contribution. This gives the first
 row. Next,
 $$
 T_n^- L=L\cdot \frac{\Psi_{n-1}}{\Phi_n},\qquad 
 R_{kl} T_n^- L= (R_{kl} L)\cdot \frac{\Psi_{n-1}}{\Phi_n}
 +L\cdot R_{kl}\Bigl( \frac{\Psi_{n-1}}{\Phi_n}\Bigr)
 $$
The second term in the right-hand side is 0 for $l<n$. For $l=n$
this term equals $R_{kn} \Psi_{n-1}/\Phi_n$. 
\hfill $\square$

\begin{lemma}
	\begin{equation}
	\Theta_{1n}= \Bigl(\sum_{J=\{j_1<j_2<\dots<j_s\}\triangleleft\{1,2,\dots n\}}	
	\zeta_{j_1j_2}\zeta_{j_2j_3}\dots \zeta_{j_{s-1}j_s}\Bigr) L\cdot \frac{\Psi_{n-1}}{\Phi_n}
	.
	\end{equation}
\end{lemma}

{\sc Proof.} First, we evaluate
$$R_I (T_n^-L)\cdot L^{-1}.$$
Each $R_I$ is a product $R_{i_1i_2}\dots R_{i_{s-1}i_s}$.
Each $R_{i_t i_{t+1}}$ is a first order differential operator without term of order zero,
therefore we can expand $R_I L$ according the Leibniz rule. 

Many summands of this expansion are zero
by a priory reasons.
Indeed, $R_{kl}\Phi_\mu$, $R_{kl}\Psi_\mu$ are nonzero only if $k\le \mu<l$.
Also, $R_{kl}R_{ab} \Phi_\mu = 0$  if $l<a$. In the case $l=a$ we have
$R_{kl}R_{lb} \Phi_\mu = -R_{kb}\Phi_\mu $.

  Therefor the  expression
$(R_I T_n^- L)/L$ is a sum of products of the following type
\begin{equation}
A[m_1,m_2]\,A[m_2,m_3] \dots A[m_{t-1},m_t] \cdot \frac{\Psi_{n-1}}{\Phi_n},
\label{eq:A[]}
\end{equation}
where $M=\{m_1,\dots, m_t\}\triangleleft I$.
 and each factor $A[m_\tau,m_{\tau+1} ]=A[i_a,i_c]$
has a form
\begin{multline*}
\frac{(q_\mu-r_{\mu+1})\,R_{i_a i_{a+1}} R_{i_{a+1} i_{a+2}}\dots R_{i_{c-1}i_c} \Phi_\mu }{\Phi_{\mu}}
\\
\text{or}\quad 
\frac{(r_{\mu+1}-q_{\mu+1})\,R_{i_a i_{a+1}} R_{i_{a+1} i_{a+2}}\dots R_{i_{c-1}i_c} \Psi_\mu }{\Psi_{\mu}}
,\end{multline*}
where $i_{c-1}\le\mu<i_c$. These expressions are equal  correspondingly
$$
\frac{(-1)^{c-a+1}(q_\mu-r_{\mu+1})\,R_{i_a i_c} \Phi_\mu }{\Phi_{\mu}}
\qquad\text{and}\qquad
\frac{(-1)^{c-a+1}(r_{\mu+1}-q_{\mu+1})\,R_{i_a i_c} \Psi_\mu }{\Psi_{\mu}}
$$
Now fix $M$ and consider the sum $\cS[I,M]$ of all summands (\ref{eq:A[]}) with fixed $M$.
It is easy to see that
$$
\cS[I,M]=(-1)^{c-a+1}\biggl[\prod_{\gamma=1}^{t-1} \Bigl(\sum_{i_a\in I: m_\gamma\le i_a<m_{\gamma+1}}
\sum_{\mu:i_a\le \mu< m_c} W_{i_a m_{\gamma+1}}^\mu\biggr]\cdot L.
$$
Next, we represent $\Theta_{1n}$ 
as
\begin{multline*}
\Theta_{1n}=\Bigl(\sum_{I\triangleleft [1,n]}
\prod_{l\in [1,n]\setminus I} (l-n+q_l-q_n)\cdot
\sum_{M\triangleleft I} \cS[I,M]\Bigr)\cdot L\cdot \frac{\Psi_{n-1}}{\Phi_n}
=\\=
\sum_{M \triangleleft [1,n]}\Bigl( \sum_{I\supset M} \prod_{l\in [1,n]\setminus I} (l-n+q_l-q_n) 
 \cdot \cS[I,M]\Bigr)\cdot L\cdot \frac{\Psi_{n-1}}{\Phi_n}.
\end{multline*}
Fix $M=\{m_1<\dots<m_t\}$. Each summand in the big brackets splits into a product of the form
$$
\prod_{\gamma=1}^{t-1} H[m_\gamma, m_{\gamma+1}],
$$
where $H[m_\gamma, m_{\gamma+1}]$ is an expression, which depend on $m_\gamma$, $m_{\gamma+1}$
and does not depend on other $m_i$. Since the coefficients also are multiplicative in the same sense,
the whole sum in the big brackets also splits in a product
$$
\prod_{\gamma=1}^{t-1} \cZ[m_\gamma, m_{\gamma_1}],
$$
where each factor $\cZ[m_\gamma, m_{\gamma+1}]$ depends only on 
$m_\gamma$, $m_{\gamma+1}$ and not on the remaining elements of a set $M$.

It remains to evaluate factors $\cZ[\cdot]$, 
$$
\cZ[k,l]=\sum_{\mu=k}^l \upsilon_{kl}^\mu W^\mu_{kl},
$$
where
\begin{multline*}
\upsilon_{kl}^\mu=
\sum_J \Bigl(-\prod_{m\in [k,\mu]\setminus J}(-1) (m-n+q_n-q_m)\cdot \prod_{m=\mu+1}^{l-1} 
(m-n+q_n-q_m)\Bigr)  
=\\=
\prod_{m=\mu+1}^{l-1} 
(m-n+q_n-q_m)
\cdot
\Bigl(-\sum_J \prod_{m\in [k,\mu]\setminus J} (-1)(m-n+q_n-q_m)
 \Bigr)        
,\end{multline*}
 the summation here is taken other all subsets $J\subset [k,\mu]$,
containing $k$. In the brackets we get
$$\prod_{m=k+1}^\mu (-1)\bigl[1-(m-n+q_n-q_m)\bigr]
$$
Thus  $\cZ_{kl}$ equals to $\zeta_{kl}$.
\hfill $\square$

\sm

Let $\zeta_{kl}$ be as above {\rm (\ref{eq:zeta})}.	Denote 
	$$
	\Theta_{mn}:=	\Bigl(\sum_{J=\{j_1<j_2<\dots<j_s\}\triangleleft\{m,2,\dots n\}}	
	\zeta_{j_1j_2}\zeta_{j_2j_3}\dots \zeta_{j_{s-1}j_s}\Bigr) L\cdot \frac{\Psi_{n-1}}{\Phi_n}
.	$$
By the previous lemma, this notation is compatible with the earlier notation $\Theta_{1n}$. Also,
$\Theta_{nn}= \frac{\Psi_{n-1}}{\Phi_n}$.

\begin{lemma}
The $\Theta_{mn}$ satisfies the following recurrence relation,
$$
\Theta_{mn}=\zeta_{m(m+1)}\Theta_{(m+1)n}+ \zeta_{m(m+2)}\Theta_{(m+2)n}+
\dots +\zeta_{mn}\Theta_{nn}.
$$		
\end{lemma}

This statement is obvious.

\sm

{\bf\punct Evaluation of $\Theta_{1n}$.\label{ss:Theta1n}}
Denote 

\begin{align*}
 s_{m\tau}^\mu&=\prod_{j=m}^{\mu-1} (q_\tau-q_j+j-n)
 \\
 \sigma^\mu_\tau&=  (q_\mu-r_{\mu+1})\prod_{i=\mu+2}^{\tau}(q_\tau-r_i-n-1+i).
\end{align*}

\begin{lemma}
\label{l:Theta-mn}
$$
\Theta_{mn} L^{-1}=-\sum_{\mu=m}^{n} s_{mn}^\mu \sigma^\mu_n \frac{R_{m\mu}\Psi_{\mu-1}}{\Phi_\mu}.
$$	
\end{lemma}

In particular, this gives an explicit expression for $\Theta_{1n}$.

\sm

{\sc Proof.}  We prove our statement by  induction. Assume that for
$\Theta_{nn}$, $\Theta_{(n-1)n}$, \dots, $\Theta_{(m+1)n}$ the formula is correct. We must
derive the equality
\begin{multline*}
\Theta_{mn}\cdot L^{-1}=\sum_{\gamma=m+1}^n \zeta_{m\gamma} \Theta_{(m+1)n}=
\\=\sum_{\gamma=m+1}^n \Bigl(\sum_{\mu=m}^{\gamma-1}w_{m\gamma}^\mu W_{m\gamma}^\mu \Bigr)
\Bigl(\sum_{\nu=\gamma}^{n} s_{\gamma n}^{\nu} \sigma^\nu\frac{R_{\gamma\nu}\Psi_{\nu-1}}{\Phi_\nu}\Bigr)
.
\end{multline*}

We have
\begin{multline*}
\upsilon_{m\gamma}^\mu s_{\gamma n}^n
=\prod_{j=m+1}^{\mu} (q_n-q_j+j-n+1) \cdot \prod_{j=\mu+1}^{\gamma-1} (q_n-q_j+j-n+1)
\prod_{j=\gamma}^{\nu-1} (q_n-q_j+j-n)
=\\=\prod_{j=m+1}^{\nu-1} (q_n-q_j+j-n+1) \cdot \prod_{j=\mu+1}^{\gamma-1} (q_n-q_j+j-n+1)
=:\xi_{m\nu}^\mu.
\end{multline*}
We stress that the expression $\xi_{m\nu}^\mu$ does not depend on $\gamma$. Our sum transforms to 
$$
\sum_{\mu,\nu: \mu<\nu}  \xi_{m\nu}^\mu \sigma^\nu_n
\biggl[\sum_{\gamma=\mu+1}^{\nu} W_{m\gamma}^{\mu} \frac{R_{\gamma\nu}\Psi_{\nu-1}}{\Phi_\nu}  \biggr]
+\frac{\Psi_{n-1}}{\Phi_n}.
$$ 
Denote by $B_{m\nu}^\mu$ the expression in the square brackets and write it explicitly:
\begin{multline}
B_{m\nu}^\mu=\sum_{\gamma=\mu+1}^{\nu}\Bigl(
\frac{(r_{\mu+1}-q_{\mu+1})R_{m\gamma}\Psi_\mu}{\Psi_{\mu}}
+\frac{(q_\mu-r_{\mu+1}) R_{m\gamma}\Phi_\mu}{\Phi_\mu}
\Bigr)\frac{R_{\gamma\nu}\Psi_{\nu-1}}{\Phi_\nu}
=\\=
(r_{\mu+1}-q_{\mu+1})\sum_{\gamma=\mu+1}^{\nu}
\frac{R_{m\gamma}\Psi_\mu}{\Psi_{\mu}} \frac{R_{\gamma\nu}\Psi_{\nu-1}}{\Phi_\nu}+
\\+
(q_\mu-r_{\mu+1})\sum_{\gamma=\mu+1}^{\nu} \frac{R_{m\gamma}\Phi_\mu}{\Psi_\mu}\frac{R_{\gamma\nu}\Psi_{\nu-1}}{\Phi_\nu}
.
\end{multline}
Next, we apply the quadratic relations
(\ref{eq:Plu-1}) and (\ref{eq:Plu-2})
and transform the last expression to
\begin{multline*}
 (r_{\mu+1}-q_{\mu+1})
 \frac{-R_{m\nu }\Psi_{\nu-1}}{\Phi_\nu}
 +(q_\mu-r_{\mu+1})\Bigl(- \frac{R_{m\nu}\Psi_{\nu-1}}{\Phi_\nu}+
 \frac{R_{m\nu} \Psi_{\mu-1}}{\Phi_\mu}\Bigr)
 =\\=
 (q_{\mu+1}-q_\mu) \frac{R_{m\nu}\Psi_{\nu-1}}{\Phi_\nu}+
 (q_\mu-r_{\mu+1}) \frac{R_{m\nu} \Psi_{\mu-1}}{\Phi_\mu}
 .
\end{multline*}
Thus,
\begin{equation}
\Theta_{mn}:=
\sum_{\mu,\nu: \mu<\nu}  \xi_{m\nu}^\mu \sigma^\nu_n
\Bigl((q_{\mu+1}-q_\mu) \frac{R_{m\nu}\Psi_{\nu-1}}{\Phi_\nu}+
 (q_\mu-r_{\mu+1}) \frac{R_{m\nu} \Psi_{\mu-1}}{\Phi_\mu}
 \Bigr)
 \label{eq:Theta-mn}
.
\end{equation}
We collect similar terms and get
$$
\Theta_{mn}:=
\sum_{\nu} \Bigl(S_1+S_2\Bigr) \frac{R_{m\nu}\Psi_{\nu-1}}{\Phi_\nu},
$$
where
\begin{align}
 S_1:&= \sum_{\mu=1}^{\nu-1} \xi_{\mu\nu}\sigma^{\nu}_n(q_{\mu+1}-q_\mu) ,
 \\
 S_2:&=\sum_{\kappa=\nu+1}^{n} \xi_{\nu\kappa} \sigma^{\kappa}_n(q_\nu-r_{\nu+1}) 
 .
\end{align}
First, we transform  $S_2$,
\begin{multline}
 S_2= (q_\nu-r_{\nu+1})\cdot \prod_{j=m+1}^{\nu} (q_n-q_j-n+j-1) 
 \times\\\times
 \sum_{\kappa=\nu+1}^{n}\Bigl( \prod_{j=\nu+1}^{\kappa-1} (q_n-q_j+j-n) (q_\kappa-r_{\kappa+1})\cdot
 \prod_{j=2}^{n-\kappa}(q_n-r_{n+2-j}-j)\Bigr).
\end{multline}

The sum in the second line is evaluated in the following lemma.

\begin{lemma}
 \label{l:sublemma-1}
 Let $n-k\ge\nu+1$. Then 
 \begin{multline}
  \sum_{\kappa=n-k+1}^{n} \prod_{j=\nu+1}^{\kappa-1} (q_n-q_j+j-n)       \cdot
  (q_\kappa-r_{\kappa+1})\cdot
 \prod_{j=2}^{n-\kappa}(q_n-r_{n+2-j}-j)=
 \\=
 \prod_{\nu+1}^{n-k} (q_n-q_j+j-n)\cdot \prod_{j=2}^k (q_n-r_{n+2-j}-k+1)
 .
 \end{multline}
 In particular,
 $$
 \sum_{\kappa=\nu+1}^{n}(\dots)= \prod_{j=2}^{n-\nu} (q_n-r_{n+2-j}-k+1)
 ,
 $$
 and
 \begin{equation}
 (q_\nu-r_{\nu+1}) \sum_{\kappa=\nu+1}^{n}(\dots)=\sigma^\nu_n.
 \end{equation}
\end{lemma}

{\sc Proof of Lemma \ref{l:sublemma-1}.}
We prove the statement by induction. Let for a given $k$ the statement hold. Then
\begin{multline*}
  \sum_{\kappa=n-k}^{n}= \sum_{\kappa=n-k+1}^{n}+\Bigl\{\text{$n-k$-th term} \Bigr\}
  =
  \\= 
  \prod_{\nu+1}^{n-k} (q_n-q_j+j-n)\cdot \prod_{j=2}^{k} (q_n-r_{n+2-j}-k+1)+\\+
  \prod_{j=\nu+1}^{n-k-1} (q_n-q_j+j-n) \cdot (q_{n-k}-r_{n-k+1})\cdot
 \prod_{j=2}^{k}(q_n-r_{n+2-j}-j)
 = \\=
 \prod_{j=\nu+1}^{n-k-1} (q_n-q_j+j-n) \prod_{j=2}^k (q_n-r_{n+2-j}-k+1)
 \cdot\Bigl( (q_n-q_{n-k}-k)+ (q_{n-k}-r_{n-k+1}) \Bigr)
\end{multline*}
and the big bracket joins to the product $\prod_{j=2}^k$.
\hfill $\square$

\medskip

Let us return to our calculation.
We must evaluate
$S_1+S_2$,
\begin{multline}
 S_1+S_2=
  \sigma^{\nu}_n\sum_{\mu=1}^{\nu-1} \xi_{\mu\nu}(q_{\mu+1}-q_\mu) +S_2
  =\\
 \sigma_{\nu} \sum_{\mu=1}^{\nu-1} \Bigl\{\prod_{j=m+1}^\mu (q_n-q_j+j-n-1) \cdot \prod_{j=\mu+1}^{\nu-1}
 (q_n-q_j+j-n)\cdot (q_\mu-q_{\mu+1})\Bigr\}
 +\\+
 \sigma^\nu_n \prod_{j=m+1}^\nu (q_n-q_j+j-n-1).
 \label{eq:last-summand}
\end{multline}
An evaluation of this sum is similar to the proof of Lemma \ref{l:sublemma-1}. We verify the following identity
by induction,
$$
 \sigma^{\nu}_n \sum_{\nu-k}^{\nu-1} \xi_{\mu\nu}(q_{\mu+1}-q_\mu) +S_2=
 \prod_{j=m+1}^{\nu-k+1}(q_n-q_j+j-1-n) \prod_{j=\nu-k+1}^{\nu-1}(q_n-q_j+j-n).
$$
Substituting $k=\nu-1$ we get the desired coefficient in (\ref{eq:Theta-mn}).

\sm

It remains to notice that formulas for the coefficient in the front of $\frac 1{\Phi_n}$  are slightly different.
In this case, the sum $S_2$ is absent, but $\zeta_{mn}$ has an extra term $\frac 1{\Phi_n}$.
Starting this place, the calculation is the same, the extra term replaces  the additional
term in (\ref{eq:last-summand}). 

\sm

This completes a proof of Lemma \ref{l:Theta-mn}.

\sm

{\bf\punct The extension of calculations.} 
Next, we must evaluate other summands in $F_{1(n+1)}L$.
Denote 
$$
\Theta_{1\tau}:= \sum_{I\triangleleft [1,\tau]}
	 \prod_{l\in [1,\tau]\setminus I} (l-\tau+q_l-q_\tau) \cdot R_I
	\Bigr) T^-_m L.
$$

\begin{lemma}
 \begin{equation}
\Theta_{1\tau}= 
	-\sum_{\mu=1}^{m} s_{1\tau}^\mu \sigma^\mu_\tau \frac{R_{1\mu}\Psi_{\mu-1}}{\Phi_\mu}
	.
	\end{equation}
\end{lemma}

{\sc Proof.} We decompose $T_\tau^-L$ as 
\begin{multline*}
T_\tau^- L=	\Bigl(\Phi_1^{-r_2+q_1}\Psi_1^{r_2-q_2}\dots \Psi_{\tau-1}^{r_{\tau}-q_{\tau}+1} \Phi_\tau^{-r_{\tau+1}+q_\tau-1}\Bigr) 
\times\\\times
\Bigl( \Psi_{m-1}^{r_{m}-q_{m}} 
\Phi_{\tau}^{-r_{\tau+1}+q_{\tau}} 
\dots \Phi_n^{-r_{n+1}+q_n}\Bigr).
\end{multline*}
Factors $R_{kl}$ of operators $R_I$ act as zero on the factors of the second bracket. Therefore this bracket can
be regarded as a constant. After this, we get the same calculations as for $\Theta_{1n}$ but $n$ is replaced by
$\tau$.
\hfill $\square$

\sm

{\bf \punct The summation.} 

\begin{lemma}
 $$
 F_{1(n+1)} L+ E_{1(n+1)} L=0.
	 $$
\end{lemma}

{\sc Proof.} 
$$
F_{1(n+1)}L=\sum_{\tau=1}^n A_\tau \sum_{\alpha=1}^\tau s_{1\tau}^\alpha \sigma_\tau^\alpha  \frac{   R_{1\alpha}\Psi_{\alpha-1}}{\Phi_\alpha}
=\sum_{\sigma=1}^n \Bigl(\sum_{\tau=\alpha}^n A_\tau s_{1\tau}^\alpha  \sigma_\tau^\alpha  \Bigr) \frac {R_{1\alpha} \Psi_{\alpha-1}}{\Phi_\alpha}.
$$
We must verify the following identities 
\begin{equation}
\sum_{\tau=\alpha}^n
A_\tau s_{1\tau}^\alpha  \sigma_\tau^\alpha =-(q_\alpha-r_{\alpha+1})
\label{eq:desired}
.\end{equation}
A summand in the left-hand side equals
\begin{multline*}
A_\tau s_{1\tau}^\alpha  \sigma_\tau^\alpha=
\frac
{\prod_{j=1}^{n-\tau+1} (q_m-r_{n+2-j}+n-\tau+1-j) }
{\prod_{j\ne \tau} (q_\tau-q_j+j-m)}
\times\\\times
(q_\alpha-r_{\alpha+1})
\prod_{i=1}^{\alpha-1} (q_m-q_i+i-m)\cdot  \prod_{k=n-m+2}^{n-\alpha} (q_m-r_{n+2-k}+n-\alpha-k+1)
.
\end{multline*}
The factor $\prod_{i=1}^{\alpha-1} (q_m-q_i+i-m)$ cancels, the factors $\prod_{j=1}^{n-\tau+1}$
and $\prod_{k=n-m+2}^{n-\alpha} $ join together and we come to 
$$
A_\tau s_{1\tau}^\alpha  \sigma_\tau^\alpha=
(q_\alpha-r_{\alpha+1})
\frac{\prod_{j=1}^{n-\sigma} (q_m-r_{n+2-j}+n-\tau+1-j}{\prod_{j\ne \tau,j>\sigma}(q_\tau-q_j+j-m)}
.
$$
Next, we write a rational function
$$
V_\alpha(x):=\frac{\prod_{j=1}^{n-\sigma} (x+r_{n+2-j}+n-j+1)}
{\prod_{j=\sigma}^n(x-q_j+j)}
$$
and evaluate the sum of its residues,
\begin{equation}
\sum_{j=\sigma}^m \mathrm{Res}_{x=q_j-j} V_\alpha(x)=- \mathrm{Res}_{x=\infty} V_\alpha(x)
\end{equation}
Multiplying this identity by 
 $(q_\alpha-r_{\alpha+1})$ we get (\ref{eq:desired}).
 \hfill $\square$
 
\sm

{\bf \punct Invariance of the space $\frV_\bfr$.}

\begin{lemma}
 The space $\frV_\bfr$ {\rm(}see {\rm(\ref{eq:frV}))}  is invariant with respect to the operator
 $F_{1(n+1)}$.
\end{lemma}

{\sc Prove.} Define elements $\ell_Z(U)\in \frV_\bfr$, where $Z$ ranges in $T_{n+1}^+$, by
$$
\ell_Z(U)=L(Z,U).
$$
First,  functions $\ell_Z(U)$ generate the space $\frV_\bfr$
(because the pairing $V_\bfp\times \frV_\bfr\to\C$ determined by $L$ is nondegenerate).
Next, 
$$
F_{1(n+1)} \ell_Z(U)= F_{1(n+1)} L(Z,U)=- E_{1(n+1)} L(Z,U)=   - \frac\partial{\partial z_{1(n+1)}}\bigl(\ell_Z(U)\bigr)
.
$$
Differentiating a family of elements of $\frV_{\bfr}$ with respect to a parameter we get elements of the same space
$\frV_{\bfr}$.
\hfill $\square$

\sm

{\bf \punct Formula for $F_{(n+1)n}$.}
Here a calculation is more-or-less the same, we omit details
(in a critical moment we use the quadratic identities (\ref{eq:Plu-1}), (\ref{eq:Plu-3})).

\section{Infinite dimensional case. A conjecture}

\COUNTERS

{\bf \punct Principal series.} Now consider two collections of complex numbers
$(p_1,\dots,p_n)$ and $(p_1^\circ,\dots,p_n^\circ)$
such that 
$$p_j-p_j^\circ\in \Z, \qquad \Re (p_j+p_{j}^\circ)=-2(j-1).$$
Consider a  representation
$\rho_{\bfp|\bfp^\circ}^n$ of the group $\GL(n,\C)$ in the space $L^2(N_{n}^+)$ 
determined by the formula
\begin{equation}
\rho_{\bfp|\bfp^\circ}(g) \phi(g)= \phi(Z^{[g]})\cdot \prod_{j=1}^n b_{jj}(Z,g)^{p_j} \ov{b_{jj}(Z,g)}^{\,p_j^\circ}
\label{eq:principal}
.\end{equation}
We get  the  unitary (nondegenerate) principal series of representations of $\GL(n,\C)$
(see \cite{GN}).
Denote by $\Lambda_n^{unitary}$ the space of all parameters $p$, $p^\circ$.

\sm

{\sc Remark.}
In formula (\ref{eq:principal}), we have complex numbers in complex powers. We understand them
in the following way:
$$
b_{jj}^{p_j}\, \ov b_{jj}^{p_j^\circ}=|b_{jj}|^{2 p_j} \,(\ov b_{jj}/b_{jj})^{p_j^\circ - p_j}.
$$
In the right hand side,  the first factor has a positive base of the power,  the second factor
has an integer exponent. Hence the product is well defined.
\hfill $\boxtimes$

\sm

{\sc Remark.} Formula (\ref{eq:principal}) makes sense if $p_j-p_j^\circ\in \Z$,
and gives a non-unitary principal series of representations.
The construction of holomorphic representation discussed above
corresponds to $p\in \Lambda_n$, $p^\circ=0$. If both 
\begin{equation}
p, \,\,p^\circ\in \Lambda_n
\label{eq:finite-dim}
\end{equation}
then our representation contains a finite-dimensional (nonholomorphic)
representation $\rho_\bfp\otimes \ov\rho_{\bfp^\circ}$, where $\ov \rho$ denotes
the complex conjugate representation.
\hfill $\boxtimes$

\sm

{\bf\punct Restriction to the smaller subgroup.} 
Consider a representation $\rho_{\bfp|\bfp^\circ}^{n+1}$ of the group $\GL(n+1,\C)$.
According \cite{Anh},
the restriction of $\rho_{\bfp|\bfp^\circ}$ to the subgroup $\GL(n,\C)$ is a multiplicity free integral
of all representations of  $\rho_{\bfp|\bfp^\circ}^n$ of unitary nondegenerate principal series
of $\GL(n,\C)$. Moreover,  the restriction has Lebesgue spectrum.
Thus the restriction $U_{\bfp|\bfp^\circ}$ can be realized in the space $L^2(N_n^+\times \Lambda_n^{unitary})$
by the formula 
$$
U_{\bfp|\bfp^\circ}(g) \psi(U,\bfq|\bfq^\circ )= \psi(Z^{[g]}, \bfq|\bfq^\circ)
\cdot \prod_{j=1}^n b_{jj}(Z,g)^{q_j} \ov{b_{jj}(Z,g)}^{\,q_j^\circ}
.$$

\sm

{\bf\punct Intertwining operator.}
Next, we write an integral operator
$$L^2(N_{n+1}^+)\to L^2(N_n^+ \times \Lambda_n^{unitary})$$
by the formula 
\begin{equation}
J \phi (U,\bfq|\bfq^\circ )=\int_{N_{n+1}^+}\phi(Z) L_{\bfp,\bfq}(Z,U)\, L_{\bfp^\circ,\bfq^\circ}(Z,U)\,
d\dot Z .
\label{eq:LL}
\end{equation}

{\bf\punct Additional operators.}
Denote the operators  (\ref{eq:formula-1}), (\ref{eq:formula-1}) by $F^\bfp_{1(n+1)}$, $F^\bfp_{(n+1)n}$.
Denote by $\ov F^{\bfp^\circ}_{1(n+1)}$, $\ov F^{\bfp^\circ}_{(n+1)n}$
the operators obtained from $F^\bfp_{1(n+1)}$, $F^\bfp_{(n+1)n}$ by replacing
$$
\frac\partial{\partial z_{kl}} \,\mapsto \frac\partial{\partial \ov z_{kl}}, \qquad p_j\mapsto p_j^\circ,\qquad
q_k\mapsto q_k^\circ.
$$
Define the operators 
\begin{align*}
F_{1(n+1)}^{\bfp|\bfp^\circ}= F^\bfp_{1(n+1)}+ \ov F^{\bfp^\circ}_{1(n+1)},
\\
F^{\bfp|\bfp^\circ}_{(n+1)n}=F^\bfp_{(n+1)n}+\ov F^{\bfp^\circ}_{(n+1)n}.
\end{align*}

\begin{conjecture}
 The operators $F_{1(n+1)}^{\bfp|\bfp^\circ}$, 
 $F^{\bfp|\bfp^\circ}_{(n+1)n}$ are images of the operators $E_{1(n+1)}$, $E_{(n+1)n}$
 under the integral transform {\rm (\ref{eq:LL})}.
\end{conjecture}

The statement seems doubtless, since we have  an analytic continuation
 of our
finite-dimensional formulas
from the set (\ref{eq:finite-dim}). However, this is not an automatic corollary
of our result. 
In particular, it is necessary to find the Plancherel  measure on $\Lambda_n^{unitary}$
(i.e. a measure making the operator $J$ unitary).

 \noindent
\tt Math.Dept., University of Vienna,
 \\
 Oskar-Morgenstern-Platz 1, 1090 Wien;
 \\
\& Institute for Theoretical and Experimental Physics (Moscow);
\\
\& Mech.Math.Dept., Moscow State University.
\\
\& Institute for Information Transmission (Moscow)
\\
e-mail: hepetuh(at) yandex.ru
\\
URL:www.mat.univie.ac.at/$\sim$neretin

\end{document}